\documentclass[a4paper, 11pt]{article}
\usepackage[mathscr]{eucal}
\usepackage{amssymb}
\usepackage{latexsym}
\usepackage{amsthm}
\usepackage{amsmath}
\usepackage[dvips]{graphicx}
\usepackage{psfrag}
\usepackage{a4wide}
\usepackage{color}
\usepackage{booktabs}
\usepackage[title,titletoc]{appendix}
\usepackage{url}

\newtheorem{theorem}{Theorem}[section]
\newtheorem{proposition}[theorem]{Proposition}
\newtheorem{lemma}[theorem]{Lemma}
\newtheorem{corollary}[theorem]{Corollary}
\theoremstyle{definition}
\newtheorem{definition}[theorem]{Definition}

\newcommand{\R}{{\mathbb R}}

\newtheorem{remark}[theorem]{Remark}

\makeatletter
\@addtoreset{equation}{section}

\makeatother
\title{Almost tight spherical $2e$-designs}
\author{Alexey Glazyrin\thanks{School of Mathematical and Statistical Sciences, The University of Texas Rio Grande Valley, Brownsville, TX 78520, USA. 
{\tt alexey.glazyrin@utrgv.edu}}\hspace{5pt} and Hiroshi Nozaki\thanks{Department of Mathematics Education, 
	Aichi University of Education, 
	1 Hirosawa, Igaya-cho, 
	Kariya, Aichi, 448-8542, 
	Japan. {\tt hnozaki@auecc.aichi-edu.ac.jp}} 
}

\begin{document}

\maketitle

\renewcommand{\thefootnote}{\fnsymbol{footnote}}


    \begin{abstract}
A spherical $2e$-design $X$ in $S^{d-1}$ satisfies the Fisher-type lower bound
\[
|X| \geq M(d,e)
:=\binom{d+e-1}{e}+\binom{d+e-2}{e-1}.
\]
A spherical $2e$-design is called \emph{tight} if equality holds. Bannai and Damerell (1979) proved that no tight spherical $2e$-design exists when $e>2$ and $d>2$. 
In this paper, we investigate the existence of spherical $2e$-designs whose cardinalities are close to this bound. 

A spherical $2e$-design is said to be \emph{almost tight} if
\[
|X|=M(d,e)+1.
\]
We prove that, for each fixed integer $e>1$, almost tight spherical $2e$-designs do not exist in sufficiently large dimensions. 
Furthermore, we prove that, for every integer $2\leq e\leq 12$, no almost tight spherical $2e$-design exists in any dimension $d>2$. 

One of the ingredients of our proofs is Mimura's duality argument. We use this argument for designs of size at least $M(d,e)+2$ to extend the 3-point semidefinite program of Bachoc and Vallentin to the case of designs and certify nonexistence of designs of certain sizes in many dimensions. To our knowledge, this is the first general method for lower bounds on the sizes of spherical $2e$-designs other than the linear programming method of Delsarte, Goethals, and Seidel (1977). Using our semidefinite program we find exact certificates of nonexistence for many sizes of spherical $2e$-designs with small $e$. In particular, our semidefinite program rules out the existence of 4-designs of size 11 in $S^2$ so, combined with the nonexistence of tight and almost tight designs, we show that the smallest 4-design in $S^2$ has size 12.

\end{abstract}

\section{Introduction}

Spherical designs were introduced by Delsarte, Goethals, and Seidel
\cite{DGS77} as a spherical analogue of combinatorial designs in the framework
of Delsarte's theory of codes and designs \cite{D73}.
A finite multiset $X$ in the unit sphere $S^{d-1}$ is called a
{\it spherical $t$-design} if
\[
\frac{1}{|S^{d-1}|}\int_{x\in S^{d-1}} f(x)\,d\mu(x)
=
\frac{1}{|X|}\sum_{x\in X} f(x)
\]
holds for every polynomial $f$ in $d$ variables of total degree at most $t$,
where $\mu$ is the standard surface measure on $S^{d-1}$ and
$|S^{d-1}|=\int_{x\in S^{d-1}}1\,d\mu(x)$.

Delsarte, Goethals, and Seidel \cite{DGS77} gave a natural lower bound for the
cardinality of a spherical $2e$-design, namely
\begin{equation}\label{eq:tight}
|X| \geq M(d,e)
:=
\binom{d+e-1}{e}
+
\binom{d+e-2}{e-1},
\end{equation}
which is called the Fisher-type lower bound.
On the other hand, spherical $t$-designs are known to exist for every $t$ and
$d$ \cite{SZ84,BRV13}.
One of the fundamental problems in the theory of spherical designs is to
determine the minimum possible cardinality of a spherical $t$-design in
$S^{d-1}$ for given $t$ and $d$.

A spherical $2e$-design is called {\it tight} if it attains the Fisher-type
lower bound.
Tight spherical $2e$-designs do not exist for any $e>2$ and $d>2$
\cite{BD79}.
Tight spherical $4$-designs in $S^{d-1}$ are known to exist for
$d=2,6,22$ \cite{DGS77}, while their complete classification remains open;
see \cite{BMV05,NV13} for nonexistence results.
A tight spherical design has the structure of a $Q$-polynomial association
scheme, and its inner products are closely related to the zeros of certain
orthogonal polynomials.
This connection plays an essential role in the nonexistence results of
Bannai and Damerell \cite{BD79}.

It is then natural to study spherical designs whose cardinalities are only slightly larger than the Fisher-type lower bound. Such designs cannot be expected, in general, to carry the structure of an association scheme.

A natural first case beyond the tight bound is that of almost tight spherical
designs.
A spherical $2e$-design $X\subset S^{d-1}$ is called {\it almost tight} if
$|X|=M(d,e)+1$.
In this paper, we prove that, for each fixed integer $e>1$, almost tight
spherical $2e$-designs do not exist in sufficiently large dimensions.
For an almost tight design, the possible inner products are restricted to the
zeros of certain explicitly determined polynomials.
Combining this restriction with an integrality condition, we obtain two
polynomial equations having a common inner product as a zero.
We then show that their resultant is nonzero in sufficiently large
dimensions, yielding a contradiction.
Furthermore, we prove that, for every integer $2\leq e\leq12$, no almost tight
spherical $2e$-design exists in any dimension $d>2$.

In the latter part of the paper, we go beyond the almost tight case and study spherical $2e$-designs whose cardinalities exceed the Fisher-type lower bound by a small number of points. Related problems have been studied by Boyvalenkov and coauthors; see, for
example, \cite{BN94,BN98,BDN99}, where restrictions on the inner products of designs with cardinalities close to the Fisher-type bound are used to derive nonexistence results. Related questions concerning the potential energy of spherical designs with cardinalities close to the Fisher-type bound were studied in~\cite{BDJ21}.

Our approach is based on a duality for spherical $2$-designs originating from Mimura \cite{M90}. This duality yields explicit restrictions on the possible inner products of a spherical $2e$-design, and we incorporate these restrictions into the semidefinite programming method of Bachoc and Vallentin \cite{BV08}. The bound~(\ref{eq:tight}) was proved by the linear programming method, and its subsequent improvements by Yudin \cite{Y97} and by Boyvalenkov and Nikova \cite{BN94, BN98} use the same method with different test polynomials. Our approach is of a different kind: it uses three-point constraints, and it yields bounds in ranges of dimensions where the linear programming method is not known to improve~(\ref{eq:tight}). We obtain a variety of new lower bounds for sizes of spherical $2e$-designs using this method (see Table~\ref{tab:bounds}). In particular, we show that the smallest 4-design in $S^2$ has size 12.

The Mimura-type duality also provides bounds on inner products in $2e$-designs whose sizes are below $2M(d,e)$. Such problems are also well studied in the literature \cite{BBD99, B02, BBKS09, BDHSS15}. We analyze these bounds and find several instances in which they are exact. Interestingly, the upper bound on inner products of 4-designs obtained this way is exact when there exist tight 5-designs but it is attained not on tight 5-designs but on designs obtained from them by latitudinal rotations. In particular, this observation solves the problem of finding the largest inner product in a 12-point 4-design on $S^2$.

This paper is organized as follows. In Section~\ref{sect:prelim} we introduce the main notation and provide necessary preliminaries about spherical designs and $s$-distance sets. In Section~\ref{sect:struct} we use mappings via spherical harmonics and the result of Mimura to determine structural properties of almost tight designs. In Section~\ref{sect:almost} we prove the main results about nonexistence of almost tight designs. In Section~\ref{sect:dual} we explain the Mimura-type duality for spherical designs of larger cardinalities and develop constraints on inner products in designs. In Section~\ref{sect:3-point} we explain how to extend the 3-point semidefinite programming method to spherical designs, describe our program in detail, and discuss nonexistence results obtained via this program. In Section~\ref{sect:constraint} we consider constraints on inner products from Section~\ref{sect:dual} and give several examples in which these bounds are attained.

\section{Preliminaries}\label{sect:prelim}

We introduce the notation and results used in this paper.

A polynomial $f \in \mathbb{R}[x_1,\ldots,x_d]$ is called \emph{harmonic} if its Laplacian is zero, that is,
\[
\Delta f=0,
\qquad
\Delta=\sum_{i=1}^d \frac{\partial^2}{\partial x_i^2}.
\]
Let $\operatorname{Harm}_k(\mathbb{R}^d)$ denote the linear space of homogeneous harmonic polynomials of total degree $k$ in $\mathbb{R}[x_1,\ldots,x_d]$. Put
$
h_k=\dim \operatorname{Harm}_k(\mathbb{R}^d).
$

Let $G_k^{(d)}(x)$ be the $d$-dimensional Gegenbauer polynomial of degree $k$, defined by the three-term recurrence relation
\[
\lambda_{k+1}G_{k+1}^{(d)}(x)
=
xG_k^{(d)}(x)
-
(1-\lambda_{k-1})G_{k-1}^{(d)}(x),
\]
where
\[
\lambda_0=0,
\qquad
\lambda_k=\frac{k}{d+2k-2},
\qquad
G_0^{(d)}(x)=1,
\qquad
G_1^{(d)}(x)=dx.
\]
With this normalization, we have
$
G_k^{(d)}(1)=h_k
$ 
\cite{DGS77}.

Let $\{\varphi_{k,1},\ldots,\varphi_{k,h_k}\}$ be an orthonormal basis of
$\operatorname{Harm}_k(\mathbb{R}^d)$ with respect to the normalized 
inner product on $S^{d-1}$,
\[
( f,g )
=
\frac{1}{|S^{d-1}|}
\int_{S^{d-1}} f(x)g(x)\,d\mu(x),
\]
where $\mu$ denotes the standard surface measure on $S^{d-1}$.
For $x,y \in \mathbb{R}^d$, let $\langle x,y \rangle$ denote the standard
inner product on $\mathbb{R}^d$. 
Using Gegenbauer polynomials and this orthonormal basis, we have the following
addition formula (see \cite{BBTZ17, DGS77}). 
\begin{theorem}[Addition formula]
For any $x,y \in S^{d-1}$, we have
\[
\sum_{i=1}^{h_k}\varphi_{k,i}(x)\varphi_{k,i}(y)
=
G_k^{(d)}(\langle x,y \rangle).
\]
\end{theorem}
The following lemma will be used later.

\begin{lemma}[{\cite[Lemma 2.5]{DGS77}}] \label{lem:constant}
Let
\[
G_i^{(d)}(x)G_j^{(d)}(x)
=
\sum_{k=0}^{i+j} q_k(i,j)G_k^{(d)}(x).
\]
Then
\[
q_0(i,j)=G_i^{(d)}(1)\delta_{i,j},
\]
where $\delta_{i,j}$ is the Kronecker delta.
\end{lemma}

We recall equivalent conditions for spherical $t$-designs.

\begin{theorem}[{\cite{BBTZ17}}] \label{thm:equivalent_design}
Let $X$ be a finite multiset in $S^{d-1}$. Then the following are equivalent.
\begin{enumerate}
    \item $X$ is a spherical $t$-design.
    \item For each $k$ with $1\leq k \leq t$, we have
    \[
    \sum_{x,y \in X} G_k^{(d)}(\langle x, y \rangle)=0.
    \]
    \item For each $k$ with $1\leq k \leq t$, we have
    \[
    \frac{1}{|X|^2}\sum_{x,y\in X}\langle x,y \rangle^k
    =
    \begin{cases}
    \displaystyle
    \frac{(k-1)!!}{d(d+2)\cdots(d+k-2)}
    &\text{if $k$ is even,}\\[6pt]
    0
    &\text{if $k$ is odd.}
    \end{cases}
    \]
\end{enumerate}
\end{theorem}

For a finite set $X\subset S^{d-1}$, define
\[
A(X)=\{\langle x,y \rangle \mid x,y \in X,\ x\ne y\}.
\]
A set $X$ is called an \emph{$s$-distance set} if $|A(X)|=s$.

The following Larman--Rogers--Seidel type theorem gives certain integrality conditions for sufficiently large $s$-distance sets.

\begin{theorem}[LRS type theorem {\cite{N11}}] \label{thm:LRS}
Let $X$ be an $s$-distance set in $S^{d-1}$ with
$A(X)=\{\alpha_1,\ldots,\alpha_s\}$. If
\[
|X| \geq
2\binom{d+s-2}{s-1}
+
2\binom{d+s-3}{s-2},
\]
then, for each $i \in \{1,\ldots,s\}$,
\[
K_i=
\prod_{j\ne i}
\frac{1-\alpha_j}{\alpha_i-\alpha_j}
\]
is an integer.
\end{theorem}

The integer $K_i$ is called the \emph{LRS ratio} with respect to $\alpha_i$.

For $\alpha_i \in A(X)$ and $x_0 \in X$, the number
\[
v_i(x_0)=|\{y \in X \mid \langle x_0,y\rangle=\alpha_i\}|
\]
is called the \emph{degree} of $\alpha_i$ with respect to $x_0$.
A finite set $X \subset S^{d-1}$ is called \emph{distance invariant} if,
for each $\alpha_i \in A(X)$, the degree $v_i(x_0)$ is independent of the
choice of $x_0 \in X$.

The following theorem is essential for determining the structure of a finite
subset of $S^{d-1}$ from the number of distances and the strength as a
spherical design.

\begin{theorem}[{\cite{DGS77}}] \label{thm:t_and_s}
Let $X$ be an $s$-distance set and a spherical $t$-design in $S^{d-1}$.
Let $A(X)=\{\alpha_1,\ldots,\alpha_s\}$ and put $\alpha_0=1$. Then the
following hold.
\begin{enumerate}
    \item If $t\geq s-1$, then $X$ is distance invariant.
    \item If $t \geq 2s-2$, then $(X,\{R_i\}_{i=0}^s)$ is a symmetric
    association scheme of class $s$, where
    \[
    R_i=\{(x,y) \in X \times X \mid \langle x,y \rangle =\alpha_i\}
    \qquad (0\leq i\leq s).
    \]
\end{enumerate}
\end{theorem}

A polynomial $F(x) \in \mathbb{R}[x]$ is called an \emph{annihilator}
for a finite set $A\subset \mathbb{R}$ if $F(a)=0$ for each $a \in A$
and $F(1)=1$.

For an $s$-distance set $X \subset S^{d-1}$, we define the annihilator
polynomial $F_X(x)$ for $A(X)$ by
\[
F_X(x)=\prod_{\alpha \in A(X)}\frac{x-\alpha}{1-\alpha}.
\]
The Gegenbauer coefficients of $F_X(x)$ determine the strength of $X$ as a
spherical design, as follows.

\begin{theorem}[{\cite[Theorem 6.5]{DGS77}}] \label{thm:strength}
Let $X$ be an $s$-distance set in $S^{d-1}$. Write the Gegenbauer expansion
of $F_X(x)$ as
\[
F_X(x)=\sum_{k=0}^s f_kG_k^{(d)}(x).
\]
If $X$ is a spherical $t$-design with $t\geq s$, then
\[
f_0=f_1=\cdots=f_{t-s}=\frac{1}{|X|}.
\]
Conversely, if
\[
f_0=f_1=\cdots=f_r=\frac{1}{|X|}
\]
and
\[
f_{r+1}>0,\ldots,f_s>0
\]
for some $r\leq s$, then $X$ is a spherical $(r+s)$-design.
\end{theorem}

From the LRS type theorem and distance invariance, we obtain the following
integrality condition for $s$-distance spherical $(s-1)$-designs, which will
be used to prove the nonexistence of almost tight spherical designs.

\begin{lemma} \label{lem:|X|c_0}
Let $X$ be an $s$-distance set in $S^{d-1}$ with
$A(X)=\{\alpha_1,\ldots,\alpha_s\}$. Suppose that $X$ is a spherical
$(s-1)$-design. Let $c_0^{(i)}$ be the coefficient of $G_0^{(d)}(x)$ in the
Gegenbauer expansion of
\[
F_i(x)=\prod_{j\ne i}\frac{x-\alpha_j}{\alpha_i-\alpha_j}.
\]
If
\[
|X| \geq
2\binom{d+s-2}{s-1}
+
2\binom{d+s-3}{s-2},
\]
then $|X|c_0^{(i)}$ is an integer for each $i \in \{1,\ldots,s\}$.
\end{lemma}

\begin{proof}
Since $X$ is an $s$-distance spherical $(s-1)$-design, Theorem
\ref{thm:t_and_s} (1) implies that $X$ is distance invariant. Let $v_i$
be the degree corresponding to $\alpha_i$.

Since $X$ is a spherical $(s-1)$-design and $\deg F_i=s-1$, we have
\[
c_0^{(i)}
=
\frac{1}{|X|^2}
\sum_{x,y\in X}F_i(\langle x,y\rangle).
\]
On the other hand, by distance invariance and the identity
$F_i(\alpha_j)=\delta_{ij}$, we have
\[
\frac{1}{|X|}
\sum_{x,y\in X}F_i(\langle x,y\rangle)
=
F_i(1)+v_i.
\]
Therefore,
\[
|X|c_0^{(i)}=F_i(1)+v_i.
\]

By the LRS type theorem, the LRS ratio
\[
F_i(1)=\prod_{j\ne i}\frac{1-\alpha_j}{\alpha_i-\alpha_j}
\]
is an integer. Since $v_i$ is also an integer, it follows that
$|X|c_0^{(i)}$ is an integer.
\end{proof}

\section{The structure of almost tight spherical $2e$-designs}\label{sect:struct}

Let
\[
H=H_e=h_1+\cdots+h_e.
\]
First, we show that mapping points of a $2e$-design to the space of spherical harmonics produces a $2$-design in dimension $H$.

Let $X$ be an arbitrary multiset in $S^{d-1}$. Let $\{\varphi_{k,i}\}_{i=1,\ldots,h_k}$ be an orthonormal basis of
$\operatorname{Harm}_k(\mathbb{R}^d)$. We consider the map
$\Phi:S^{d-1}\rightarrow S^{H-1}$ defined by
\[
\Phi(x)=\frac{1}{\sqrt{H}}
(\varphi_{1,1}(x),\ldots,\varphi_{1,h_1}(x),
\varphi_{2,1}(x),\ldots,\varphi_{2,h_2}(x),\ldots,
\varphi_{e,1}(x),\ldots,\varphi_{e,h_e}(x)).
\]
This map is well-defined, since, by the addition formula, we have
\[
\langle \Phi(x),\Phi(x)\rangle
=
\frac{1}{H}\sum_{k=1}^e G_k^{(d)}(1)
=
1.
\]
Moreover, by the addition formula, the inner product
$\langle \Phi(x),\Phi(y)\rangle$ depends only on $\langle x,y\rangle$:
\begin{equation}\label{eq:inn_phi}
\langle \Phi(x),\Phi(y)\rangle
=
\frac{1}{H}
\sum_{k=1}^e \sum_{i=1}^{h_k}
\varphi_{k,i}(x)\varphi_{k,i}(y)
=
\frac{1}{H}\sum_{k=1}^e G_k^{(d)}(\langle x,y\rangle).    
\end{equation}
In particular, the configuration $\Phi(X)$ is independent of the choice of
the orthonormal bases $\{\varphi_{k,i}\}$ up to an orthogonal transformation.
Note also that the map $\Phi$ is injective. Indeed, the degree-one harmonic
part is contained in the definition of $\Phi$, and it separates points on
$S^{d-1}$.

\begin{theorem}\label{thm:harmmap}
If $X$ is a $2e$-design in $S^{d-1}$, then $\Phi(X)$ is a $2$-design in $S^{H-1}$.
\end{theorem}

\begin{proof}
If $X\subset S^{d-1}$ is a spherical $2e$-design, then
\[
\sum_{x,y\in X}\langle \Phi(x),\Phi(y)\rangle
=
\frac{1}{H}\sum_{k=1}^e \sum_{x,y\in X}
G_k^{(d)}(\langle x,y\rangle)
=
0.
\]
Furthermore,
\begin{align*}
\sum_{x,y\in X}\langle \Phi(x),\Phi(y)\rangle^2
&=
\frac{1}{H^2}
\sum_{x,y\in X}
\left(\sum_{k=1}^e G_k^{(d)}(\langle x,y\rangle)\right)^2  \\
&=
\frac{|X|^2}{H^2}
\sum_{k=1}^e G_k^{(d)}(1) \\
&=
\frac{|X|^2}{H}.
\end{align*}
Here we used Lemma~\ref{lem:constant} and the fact that $X$ is a spherical
$2e$-design; indeed, the coefficient of $G_0^{(d)}(x)$ in
\[
\left(\sum_{k=1}^e G_k^{(d)}(x)\right)^2
\]
is equal to $\sum_{k=1}^e G_k^{(d)}(1)=H$. 
Therefore, by Theorem~\ref{thm:equivalent_design}, $\Phi(X)$ is a spherical
$2$-design in $S^{H-1}$.
\end{proof}

This result allows us to determine the set $A(X)$ of inner products for an almost tight spherical $2e$-design $X$ in $S^{d-1}$. Note that in this case $|X|=H+2$. By Theorem~\ref{thm:harmmap}, $\Phi(X)$ is an $(H+2)$-point spherical $2$-design in $S^{H-1}$. $2$-designs of cardinality $H+2$ in dimension $H$ are already described in the literature \cite{M90,S93}. In particular, an $(H+2)$-point spherical $2$-design in $S^{H-1}$ exists only when $H$ is even. In this case, such a design is the union of two mutually orthogonal regular simplices $\Lambda_1$ and $\Lambda_2$, each consisting of $H/2+1$ points and spanning a subspace of dimension $H/2$.

We call the decomposition $X=X_1\sqcup X_2$, where
$X_i=\Phi^{-1}(\Lambda_i)$ for $i=1,2$, the \emph{simplex preimage
decomposition} of $X$.
In this decomposition, we have
\[
\langle \Phi(x), \Phi(y) \rangle =
\begin{cases}
    1 &\text{if $x=y$},\\
    -2/H &\text{if $x,y \in X_i$ for some $i\in\{1,2\}$ and $x\ne y$},\\
    0 &\text{if $x\in X_i$ and $y\in X_j$ with $i\ne j$}.
\end{cases}
\]

Let
\[
R_e(x)=\sum_{k=0}^e G_k^{(d)}(x).
\]
For distinct points $x,y \in X_i$, the inner product $\langle x,y\rangle$
is a root of
\[
P_e(x)=R_e(x)+1
=
\sum_{k=0}^eG_k^{(d)}(x)+1,
\]
by \eqref{eq:inn_phi} and $\langle \Phi(x),\Phi(y)\rangle=-2/H$.
Similarly, for $x \in X_1$ and $y \in X_2$, the inner product
$\langle x,y\rangle$ is a root of
\[
Q_e(x)=R_e(x)-1
=
\sum_{k=1}^eG_k^{(d)}(x).
\]
\begin{lemma} \label{lem:distance_invariant}
Let $X$ be an almost tight spherical $2e$-design in $S^{d-1}$. Then the
inner product between any two distinct points of $X$ is a zero of either
$P_e(x)$ or $Q_e(x)$. Moreover, $X$ is distance invariant.
\end{lemma}

\begin{proof}
By the above argument, the first statement is clear. Hence $X$ is an
$s$-distance set for some $s\leq 2e$. Since $X$ is a spherical
$2e$-design, we have $2e\geq s-1$. Therefore, by Theorem
\ref{thm:t_and_s} (1), $X$ is distance invariant.
\end{proof}

\begin{remark}
Since the possible values of the inner products are discrete, an almost
tight spherical design $X$ is rigid in the sense of Bannai \cite{B87}.
\end{remark}

\begin{lemma} \label{lem:Xi_distance}
Let $X$ be an almost tight spherical $2e$-design in $S^{d-1}$ with
$e\geq 2$. Let $X=X_1\sqcup X_2$ be the simplex preimage decomposition.
Then the inner product between any two distinct points of $X_i$ is a zero
of $P_e(x)$ for $i=1,2$. Moreover, if $d\geq e+1$, then $X_i$ is an
$e$-distance set for $i=1,2$.
\end{lemma}

\begin{proof}
By the preceding discussion, the first statement follows, and hence
$|A(X_i)|\leq e$. We prove $|A(X_i)|\geq e$.

If $Y\subset S^{d-1}$ is an $s$-distance set, then
\[
|Y|\leq N(d,s):=
\binom{d+s-1}{s}+\binom{d+s-2}{s-1}
=
\frac{d+2s-1}{s}\binom{d+s-2}{s-1},
\]
as proved in \cite{DGS77}.

Since
\[
|X_i|=\frac{1}{2}
\left\{
\binom{d+e-1}{e}+\binom{d+e-2}{e-1}+1
\right\}
=
\frac{N(d,e)+1}{2},
\]
it is enough to show that
\[
\frac{N(d,e)+1}{2}>N(d,e-1).
\]
Indeed, from $d\geq e+1$ and $e\geq 2$, we have
\begin{align*}
N(d,e)+1-2N(d,e-1)
&=
\left(
\frac{(d-1)(d+e-2)}{e(e-1)}-2
\right)
\binom{d+e-3}{e-2}
+1\\
&\geq
\left(
\frac{e(2e-1)}{e(e-1)}-2
\right)
\binom{d+e-3}{e-2}
+1\\
&=
\frac{1}{e-1}\binom{d+e-3}{e-2}
+1\\
&>0.
\end{align*}
Therefore, $|X_i|>N(d,e-1)$, and hence $X_i$ cannot be an
$s$-distance set with $s\leq e-1$. Since $|A(X_i)|\leq e$, it follows
that $|A(X_i)|=e$. Thus $X_i$ is an $e$-distance set.
\end{proof}

\begin{lemma} \label{lem:Xi_design}
Let $X$ be an almost tight spherical $2e$-design in $S^{d-1}$ with $e\geq 2$.
Let $X=X_1\sqcup X_2$ be the simplex preimage decomposition. If $d\geq e+1$, then
$X_i$ is a spherical $e$-design but not a spherical $(e+1)$-design for
$i=1,2$.
\end{lemma}

\begin{proof}
From $d\geq e+1$ and Lemma~\ref{lem:Xi_distance}, $X_i$ is an
$e$-distance set. The annihilator polynomial $F_{X_i}(x)$ is written as
\begin{align*}
F_{X_i}(x)
&=\frac{P_e(x)}{P_e(1)}\\
&=\frac{1}{H+2}
\left(2G_0^{(d)}(x)+G_1^{(d)}(x)+\cdots+G_e^{(d)}(x)\right).
\end{align*}
Since $|X_i|=(H+2)/2$, the coefficient of $G_0^{(d)}(x)$ in
$F_{X_i}(x)$ is
\[
\frac{2}{H+2}=\frac{1}{|X_i|},
\]
and the other Gegenbauer coefficients are positive. Therefore, by
Theorem~\ref{thm:strength}, $X_i$ is a spherical $e$-design.

On the other hand, the coefficient of $G_1^{(d)}(x)$ in $F_{X_i}(x)$ is
\[
\frac{1}{H+2}\ne \frac{1}{|X_i|}.
\]
Thus, again by Theorem~\ref{thm:strength}, $X_i$ is not a spherical
$(e+1)$-design.
\end{proof}

\begin{lemma}
Let $X$ be an almost tight spherical $2e$-design in $S^{d-1}$ with $e\geq 2$.
Let $X=X_1\sqcup X_2$ be the simplex preimage decomposition. If $d\geq e+1$, then
$X_i$ is distance invariant for $i=1,2$.
\end{lemma}

\begin{proof}
By Lemmas~\ref{lem:Xi_distance} and~\ref{lem:Xi_design}, $X_i$ is an
$e$-distance set and a spherical $e$-design. Therefore, Theorem
\ref{thm:t_and_s} (1) implies that $X_i$ is distance invariant.
\end{proof}

\section{Almost tight spherical $2e$-designs in $S^{d-1}$}\label{sect:almost}

\subsection{The cases $d=2$ and $e=1$}

In the case $d=2$, an almost tight spherical $2e$-design in $S^1$ has
$2e+2$ points. By Hong's classification of spherical $t$-designs in
$\mathbb{R}^2$ \cite{H82}, it must be a regular $(2e+2)$-gon.

In the case $e=1$, an almost tight spherical $2$-design in $S^{d-1}$ has
$d+2$ points. If $d$ is odd, then no almost tight spherical $2$-design in
$S^{d-1}$ exists \cite{M90}. If $d$ is even, then any almost tight
spherical $2$-design in $S^{d-1}$ is the union of two mutually orthogonal
regular simplices $\Lambda_1$ and $\Lambda_2$, each consisting of
$d/2+1$ points and spanning a subspace of dimension $d/2$ \cite{M90,S93}. 

\subsection{Nonexistence of almost tight spherical $4$-designs}

Gorbachev~\cite{Gor26} recently proved that there is no $10$-point
spherical $4$-design on $S^2$. The following theorem extends this
nonexistence result to almost tight spherical $4$-designs in every
dimension $d>2$.

\begin{theorem} \label{thm:4-design}
An almost tight spherical $4$-design in $S^{d-1}$ does not exist for any
$d>2$.
\end{theorem}

\begin{proof}
Let $X$ be an almost tight spherical $4$-design in $S^{d-1}$ with $d>2$.
Then
\[
|X|=\binom{d+2}{2}.
\]
Let $X=X_1\sqcup X_2$ be the simplex preimage decomposition of $X$.
Since $e=2$ and $d>2$, the condition $d\geq e+1$ in Lemmas
\ref{lem:Xi_distance} and~\ref{lem:Xi_design} is satisfied. Hence $X_i$ is
a $2$-distance set and a spherical $2$-design for $i=1,2$. Therefore,
by Theorem~\ref{thm:t_and_s} (2), $X_i$ has the structure of a strongly
regular graph.

The inner products between distinct points of $X_i$ are the zeros of
$P_2(x)$. We have
\[
P_2(x)
=
2+dx+\frac{d(d+2)}{2}\left(x^2-\frac{1}{d}\right)
=
\frac{d(d+2)}{2}
\left(
x^2+\frac{2}{d+2}x-\frac{d-2}{d(d+2)}
\right).
\]
Thus the roots of $P_2(x)$ are
\[
\alpha_1=
\frac{-d+\sqrt{d(d^2+d-4)}}{d(d+2)},
\qquad
\alpha_2=
\frac{-d-\sqrt{d(d^2+d-4)}}{d(d+2)}.
\]

Since $X_i$ has the structure of a strongly regular graph, the inner
products are rational unless the corresponding graph is a conference graph.
If the graph of $X_i$ is a conference graph, then $|X_i|=2d+1$. However,
\[
|X_i|=\frac{1}{2}\binom{d+2}{2}
=
\frac{(d+1)(d+2)}{4},
\]
which is never equal to $2d+1$ for any integer $d$. Therefore, the inner
products $\alpha_1$ and $\alpha_2$ must be rational.

For $3\leq d\leq 9$ with $d\ne 4$, the numbers $\alpha_1$ and $\alpha_2$
are irrational, which is a contradiction. For $d=4$, we have $|X|=15$,
and hence $|X_i|=15/2$, which is impossible.

Suppose $d\geq 10$. Then
\[
|X_1|=\frac{1}{2}\binom{d+2}{2}>2d+2.
\]
By Theorem~\ref{thm:LRS},
\[
K_1=\frac{1-\alpha_2}{\alpha_1-\alpha_2}
\]
is an integer. Put $K_1=m$. Then
\[
\frac{d(d+3)}{\sqrt{d(d^2+d-4)}}=2m-1\in\mathbb{Z}.
\]
Hence its square is also an integer. On the other hand,
\[
\left(
\frac{d(d+3)}{\sqrt{d(d^2+d-4)}}
\right)^2
=
d+5+\frac{8d+20}{d^2+d-4}.
\]
For $d\geq 10$, we have
\[
0<\frac{8d+20}{d^2+d-4}<1.
\]
Thus the above square is not an integer, a contradiction. This contradiction shows
that no almost tight spherical $4$-design exists for $d\geq 10$.

Therefore, no almost tight spherical $4$-design exists in $S^{d-1}$ for
any $d>2$.
\end{proof}

\begin{remark}
We can also prove Theorem~\ref{thm:4-design} by determining the solutions in positive integers $d,n$ of the elliptic curve equation
\[
n^2=d(d^2+d-4),
\]
namely
\[
(d,n)=(2,2),(4,8),(18,78).
\]
For $d=18$, we have
\[
\alpha_1=\frac{1}{6},
\qquad
\alpha_2=-\frac{4}{15},
\qquad
|X_i|=95.
\]
From the $1$-design property of $X_i$, we have
\[
1+\alpha_1k_1+\alpha_2(|X_i|-1-k_1)=0,
\]
which implies
\[
k_1=\frac{722}{13},
\]
a contradiction.
\end{remark}
\subsection{Nonexistence of almost tight spherical $2e$-designs for $e>2$}

Let $X$ be an almost tight spherical $2e$-design in $S^{d-1}$ with $e>2$,  
and let $X=X_1 \sqcup X_2$ be the simplex preimage decomposition.
We suppose $d\geq 3e+1$ so that
\[
|X_1|\geq
2\binom{d+e-2}{e-1}
+
2\binom{d+e-3}{e-2}.
\]
Indeed,
\[
2\left(
|X_1|
-
2\binom{d+e-2}{e-1}
-
2\binom{d+e-3}{e-2}
\right)
=
\left(
\frac{(d-2e-1)(d+e-2)}{e(e-1)}
-4
\right)
\binom{d+e-3}{e-2}
+1.
\]
If $d\geq 3e+1$, then the right-hand side is at least
\[
\frac{3}{e-1}\binom{d+e-3}{e-2}+1>0.
\]
On the other hand, if $3\leq d\leq 3e$, then the right-hand side is negative
for $e>2$.

We first outline the idea of the proof here.
By Lemma~\ref{lem:Xi_distance}, $X_1$ is an $e$-distance set and
$P_e(\alpha)=0$ for each $\alpha\in A(X_1)$. 
For each $\alpha \in A(X_1)$, let
\[
F_\alpha(x)=
\prod_{\substack{\beta\in A(X_1)\\ \beta \ne \alpha}}
\frac{x-\beta}{\alpha-\beta},
\]
and write its Gegenbauer expansion as
\[
F_\alpha(x)=\sum_{j=0}^{e-1}b_jG_j(x).
\]
By Lemma~\ref{lem:|X|c_0} and the assumption $d\geq 3e+1$, the number
$|X_1|b_0$ is an integer. Moreover, we can show that
$b_0\in \mathbb{Q}(\alpha)$. Therefore, the integrality condition on
$|X_1|b_0$ gives a polynomial relation $g(\alpha)=0$ for some
$g(x)\in \mathbb{Q}[x]$. 
We then prove that $\operatorname{Res}_x(P_e(x),g(x))\ne 0$, which
contradicts the fact that $\alpha$ is a common root of $P_e(x)$ and
$g(x)$. This proves the nonexistence.

We give explicit expressions for the coefficients of
\[
R_n(x)=\sum_{k=0}^n G_k(x).
\]
From \cite[Equation (2.7)]{BD79}, we have
\[
R_n(x)=C_n^{(d/2)}(x)+C_{n-1}^{(d/2)}(x),
\]
where
\begin{align}
C_n^{(d/2)}(x)
&=\sum_{j=0}^{\lfloor n/2\rfloor}
(-1)^j
\frac{(d/2)_{n-j}}{j!(n-2j)!}(2x)^{n-2j} \nonumber \\
&=\sum_{j=0}^{\lfloor n/2\rfloor}
(-1)^j
\frac{d(d+2)(d+4)\cdots (d+2n-2j-2)}
{j!(n-2j)!2^j}
x^{n-2j} \nonumber \\
&=
d(d+2)\cdots (d+2n-2\lfloor n/2 \rfloor -2) \nonumber \\
&\qquad \times
\sum_{j=0}^{\lfloor n/2\rfloor}
(-1)^j
\frac{
(d+2n-2\lfloor n/2\rfloor)
(d+2n-2\lfloor n/2\rfloor+2)
\cdots
(d+2n-2j-2)
}
{j!(n-2j)!2^j}
x^{n-2j},
\label{eq:c_i}
\end{align}
where empty products are interpreted as $1$.

\begin{theorem}\label{thm:nonexistence_almost_tight}
For any integer $e>2$, there exists an integer $d_0\geq 3e+1$ such that,
for any $d\geq d_0$, an almost tight spherical $2e$-design in $S^{d-1}$
does not exist.
\end{theorem}
\begin{proof}
  We use the notation introduced above. Suppose that $d\geq 3e+1$.
By Lemmas~\ref{lem:Xi_distance} and~\ref{lem:Xi_design}, $X_1$ is an
$e$-distance set and a spherical $e$-design. Moreover, the assumption
$d\geq 3e+1$ implies
\[
|X_1|\geq
2\binom{d+e-2}{e-1}
+
2\binom{d+e-3}{e-2}.
\]
Hence Lemma~\ref{lem:|X|c_0} can be applied to $X_1$. 

By Lemma~\ref{lem:Xi_distance}, $X_1$ is an $e$-distance set, and each
element of $A(X_1)$ is a zero of $P(x)=R_e(x)+1$. Since $P(x)$ has degree
$e$, the elements of $A(X_1)$ are precisely the zeros of $P(x)$. In
particular, $P'(\alpha)=R'_e(\alpha)\ne 0$ for each $\alpha\in A(X_1)$.

For $\alpha \in A(X_1)$, the polynomial $F_\alpha(x)$ can be expressed as
\[
F_\alpha(x)
=
\frac{R_e(x)+1}{x-\alpha}\cdot\frac{1}{R_e'(\alpha)}
=
\tilde{F}_\alpha(x)\cdot \frac{1}{R_e'(\alpha)}.
\]
For a polynomial $f$, let $(f)_0$ denote the coefficient of $G_0$ in its
Gegenbauer expansion. By Lemma~\ref{lem:|X|c_0},
\[
|X_1|b_0
=
|X_1|(F_\alpha)_0
=
|X_1|(\tilde{F}_\alpha)_0\cdot \frac{1}{R_e'(\alpha)}
=
m \in \mathbb{Z},
\]
and hence
\begin{equation}\label{eq:key}
|X_1|(\tilde{F}_\alpha)_0=mR_e'(\alpha).
\end{equation}

We prove that $(\tilde{F}_\alpha)_0$ can be written as a polynomial in
$\alpha$ with coefficients in $\mathbb{Q}[d]$. Thus there exists a
polynomial $h(x)\in\mathbb{Q}[d][x]$ such that
\[
h(\alpha)=(\tilde{F}_\alpha)_0.
\]
It follows from \eqref{eq:key} that the polynomial
\[
g(x)=|X_1|h(x)-mR'_e(x)\in \mathbb{Q}[d][x]
\]
satisfies $g(\alpha)=0$. 
It remains to prove that the resultant of $P(x)=R_e(x)+1$ and $g(x)$ is
nonzero for sufficiently large $d$.

We treat the even and odd cases of $e$ simultaneously. Let
\[
e=2t+\epsilon,\qquad \epsilon\in\{0,1\},
\]
where $t\geq 1$, and assume $e\geq 3$. Put
\[
D_t=d(d+2)\cdots(d+2t-2).
\]
Define
\[
c_{e-2j}
=
(-1)^j
\frac{(d+2t)(d+2t+2)\cdots(d+2e-2j-2)}
{j!(e-2j)!2^j}
\qquad (0\leq j\leq t)
\]
and
\[
c_{e-1-2j}
=
(-1)^j
\frac{(d+2t)(d+2t+2)\cdots(d+2e-2j-4)}
{j!(e-1-2j)!2^j}
\qquad (0\leq j\leq t+\epsilon-1).
\]
Then
\[
C_e^{(d/2)}(x)
=
D_t\sum_{j=0}^{t}c_{e-2j}x^{e-2j}
\]
and
\[
C_{e-1}^{(d/2)}(x)
=
D_t\sum_{j=0}^{t+\epsilon-1}c_{e-1-2j}x^{e-1-2j}.
\]
Thus we can write
\[
R_e(x)=D_t\sum_{r=0}^{e}c_rx^r.
\]

Let $q_i(\alpha)\in\mathbb Q[d][\alpha]$ be the coefficient of $x^i$ in
$\tilde{F}_\alpha(x)$:
\[
\tilde{F}_\alpha(x)=\sum_{i=0}^{e-1}q_i(\alpha)x^i.
\]
Since $R_e(\alpha)+1=0$, dividing $R_e(x)+1$ by $x-\alpha$ gives
\[
q_i(\alpha)
=
D_t\sum_{j=0}^{e-1-i}c_{e-j}\alpha^{e-1-i-j}.
\]
Here the additional constant term $1$ in $R_e(x)+1$ does not contribute to
the coefficients $q_i(\alpha)$ of the quotient; it only affects the
remainder, which is zero because $R_e(\alpha)+1=0$.

We use $(x^{2i})_0=(2i-1)!!/(d(d+2)\cdots(d+2i-2))$ and
$(x^{2i+1})_0=0$. Therefore, we obtain
\begin{align*}
(\tilde{F}_\alpha)_0
&=\sum_{i=0}^{e-1}q_i(\alpha)(x^i)_0\\
&=\sum_{i=0}^{t+\epsilon-1}q_{2i}(\alpha)(x^{2i})_0\\
&=\sum_{i=0}^{t+\epsilon-1}q_{2i}(\alpha)
\frac{(2i-1)!!}{d(d+2)\cdots(d+2i-2)}\\
&=\sum_{i=0}^{t+\epsilon-1}
(2i-1)!!
(d+2i)(d+2i+2)\cdots(d+2t-2)
\sum_{j=0}^{e-1-2i}
c_{e-j}\alpha^{e-1-2i-j},
\end{align*}
where $(2i-1)!!=1$ for $i=0$, and empty products are interpreted as $1$.

For $0\leq i\leq t+\epsilon-1$, let
\[
\ell_{2i}
=
(2i-1)!!
(d+2i)(d+2i+2)\cdots(d+2t-2).
\]
Then, for $0\leq k\leq e-1$, the coefficient of $\alpha^k$ in
$(\tilde{F}_\alpha)_0$ is
\[
\tilde{h}_k(d)
=
\sum_{\substack{
0\leq i\leq t+\epsilon-1\\
0\leq j\leq e-1-2i\\
e-1-2i-j=k
}}
\ell_{2i}c_{e-j}.
\]
Equivalently, since $2i=e-1-j-k$, we have
\[
\tilde{h}_k(d)
=
\sum_{\substack{
0\leq j\leq e-1-k\\
j\equiv e-1-k \pmod 2
}}
\ell_{e-1-j-k}c_{e-j}.
\]
Thus
\[
h(x)=\sum_{k=0}^{e-1}\tilde{h}_k(d)x^k.
\]
By construction, we have $h(\alpha)=(\tilde{F}_\alpha)_0$.

We now compute $\operatorname{Res}_x(P(x),g(x))$. 
We record the behavior modulo $d$ of several coefficients of $P(x)$ and
$g(x)$. Throughout the following computation, $d$ and $m$ are treated as indeterminates, and congruences modulo $d$ are taken in the polynomial ring $\mathbb{Q}[d,m]$. 

Write
\[
P(x)=\sum_{i=0}^{e}p_i x^i.
\]
Since
\[
P(x)=C_e^{(d/2)}(x)+C_{e-1}^{(d/2)}(x)+1
\equiv 1\pmod d,
\]
all coefficients of $P(x)$ except the constant term are congruent to $0$
modulo $d$, and
\[
p_0\equiv 1\pmod d.
\]
Moreover,
\[
\frac{p_e}{d}
=
(d+2)(d+4)\cdots(d+2t-2)c_e
\not\equiv 0\pmod d
\]
and
\[
\frac{p_{e-1}}{d}
=
(d+2)(d+4)\cdots(d+2t-2)c_{e-1}
\not\equiv 0\pmod d.
\]

Write
\[
g(x)=|X_1|h(x)-mR'_e(x)=\sum_{i=0}^{e-1}g_i(d,m)x^i.
\]
Since $R'_e(x)\equiv 0\pmod d$ and $|X_1|\not\equiv 0\pmod d$, we have
\[
g_i(d,m)\equiv |X_1|\tilde{h}_i(d)\pmod d.
\]

For $k=e-1$, the coefficient of $\alpha^{e-1}$ in
$(\tilde{F}_\alpha)_0$ is
\[
\tilde{h}_{e-1}(d)
=
\ell_0c_e
=
D_t c_e.
\]
Hence
\[
\tilde{h}_{e-1}(d)\equiv 0\pmod d,
\qquad
\frac{\tilde{h}_{e-1}(d)}{d}\not\equiv 0\pmod d.
\]

For $k=e-2$, the coefficient of $\alpha^{e-2}$ in
$(\tilde{F}_\alpha)_0$ is
\[
\tilde{h}_{e-2}(d)
=
\ell_0c_{e-1}
=
D_t c_{e-1}.
\]
Hence
\[
\tilde{h}_{e-2}(d)\equiv 0\pmod d,
\qquad
\frac{\tilde{h}_{e-2}(d)}{d}\not\equiv 0\pmod d.
\]

For $k=e-3$, the coefficient of $\alpha^{e-3}$ in
$(\tilde{F}_\alpha)_0$ is
\[
\tilde{h}_{e-3}(d)
=
\ell_2c_e+\ell_0c_{e-2}.
\]
Since $\ell_0$ is divisible by $d$, we have
\[
\tilde{h}_{e-3}(d)\equiv \ell_2c_e\pmod d.
\]
Moreover,
\[
\ell_2c_e
=
(d+2)(d+4)\cdots(d+2t-2)c_e.
\]
Therefore,
\[
\tilde{h}_{e-3}(d)\not\equiv 0\pmod d.
\]

Now assume $e\geq 4$. We compute the resultant
$\operatorname{Res}_x(P(x),g(x))$ modulo $d$ after dividing by $d^2$. Put
\[
u_1(d,m)=\frac{g_{e-1}(d,m)}{d},
\qquad
u_2(d,m)=\frac{g_{e-2}(d,m)}{d}.
\]
The possible dependence on $m$ appears only in $u_1(d,m)$ and
$u_2(d,m)$; all the other entries in the determinant are independent of
$m$ modulo $d$.

Then, by expanding the Sylvester determinant along the entries congruent
to $1$ in the shifted rows coming from $P(x)$, we obtain
\begin{align*}
\frac{\operatorname{Res}_x(P(x),g(x))}{d^2}
&\equiv
\begin{vmatrix}
p_{e}/d & p_{e-1}/d & 0 & \cdots & 0 & 1 & 0 & \cdots & 0 \\
0 & p_{e}/d & 0 & \cdots & 0 & 0 & 1 & \cdots & 0 \\
 &  &  & \ddots &  &  &  & \ddots &  \\
0 & \cdots & 0 & 0 & 0 & 0 & 0 & \cdots & 1 \\
u_1 & u_2 & |X_1|\tilde{h}_{e-3} & \cdots & |X_1|\tilde{h}_0 & 0 & 0 & \cdots & 0 \\
0 & u_1 & 0 & |X_1|\tilde{h}_{e-3} & \cdots & |X_1|\tilde{h}_0 & 0 & \cdots & 0 \\
0 & 0 & 0 & 0 & |X_1|\tilde{h}_{e-3} & \cdots & |X_1|\tilde{h}_0 & \cdots & 0 \\
 &  &  &  & \ddots & \ddots & \ddots & \ddots &  \\
0 & \cdots & 0 & 0 & 0 & 0 & |X_1|\tilde{h}_{e-3} & \cdots & |X_1|\tilde{h}_0
\end{vmatrix}
\\
&\equiv
\begin{vmatrix}
p_e/d & p_{e-1}/d & 0 & 0 & \cdots & 0 & 1 & 0\\
0 & p_e/d & 0 & 0 & \cdots & 0 & 0 & 1\\
u_1 & u_2 & |X_1|\tilde{h}_{e-3} & |X_1|\tilde{h}_{e-4} & \cdots & |X_1|\tilde{h}_0 & 0 & 0\\
0 & u_1 & 0 & |X_1|\tilde{h}_{e-3} & |X_1|\tilde{h}_{e-4} & \cdots & |X_1|\tilde{h}_0 & 0\\
0 & 0 & 0 & 0 & |X_1|\tilde{h}_{e-3} & |X_1|\tilde{h}_{e-4} & \cdots & |X_1|\tilde{h}_0\\
\vdots &  &  &  & \ddots & \ddots & \ddots & \vdots\\
0 & \cdots & 0 & 0 & 0 & 0 & |X_1|\tilde{h}_{e-3} & |X_1|\tilde{h}_{e-4}\\
0 & \cdots & 0 & 0 & 0 & 0 & 0 & |X_1|\tilde{h}_{e-3}
\end{vmatrix}
\\
&=
\left(|X_1|\tilde{h}_{e-3}\right)^{e-2}
\begin{vmatrix}
p_e/d & p_{e-1}/d & 0 & 0\\
0 & p_e/d & 0 & 0\\
u_1 & u_2 & |X_1|\tilde{h}_{e-3} & |X_1|\tilde{h}_{e-4}\\
0 & u_1 & 0 & |X_1|\tilde{h}_{e-3}
\end{vmatrix}
\\
&=
|X_1|^e
\left(\frac{p_e}{d}\right)^2
\tilde{h}_{e-3}^{\,e}
\pmod d.
\end{align*}
Thus the reduction of
\[
\frac{\operatorname{Res}_x(P(x),g(x))}{d^2}
\]
modulo $d$ is independent of $m$. Since
\[
|X_1|\not\equiv 0\pmod d,\qquad
\frac{p_e}{d}\not\equiv 0\pmod d,\qquad
\tilde{h}_{e-3}\not\equiv 0\pmod d,
\]
we have
\[
\frac{\operatorname{Res}_x(P(x),g(x))}{d^2}
\not\equiv 0\pmod d.
\]

It remains to consider the case $e=3$, which can be checked directly by
the same argument. Put
\[
u_1(d,m)=\frac{g_2(d,m)}{d},
\qquad
u_2(d,m)=\frac{g_1(d,m)}{d}.
\]
Again, the possible dependence on $m$ appears only in $u_1(d,m)$ and
$u_2(d,m)$. Then
\[
\frac{\operatorname{Res}_x(P(x),g(x))}{d^2}
\equiv
\begin{vmatrix}
p_3/d & p_2/d & 0 & 1 & 0\\
0 & p_3/d & 0 & 0 & 1\\
u_1 & u_2 & |X_1|\tilde{h}_0 & 0 & 0\\
0 & u_1 & 0 & |X_1|\tilde{h}_0 & 0\\
0 & 0 & 0 & 0 & |X_1|\tilde{h}_0
\end{vmatrix}
\pmod d.
\]
Expanding along the last row and then along the second row, and evaluating
the remaining $3\times 3$ determinant, we obtain
\[
\frac{\operatorname{Res}_x(P(x),g(x))}{d^2}
\equiv
|X_1|^3
\left(\frac{p_3}{d}\right)^2
\tilde{h}_0^3
\pmod d.
\]
Since
\[
|X_1|\not\equiv 0\pmod d,\qquad
\frac{p_3}{d}\not\equiv 0\pmod d,\qquad
\tilde{h}_0\not\equiv 0\pmod d,
\]
we have
\[
\frac{\operatorname{Res}_x(P(x),g(x))}{d^2}
\not\equiv 0\pmod d.
\]

Therefore, we can write
\[
\operatorname{Res}_x(P(x),g(x))=d^2H(d,m)
\]
for some $H(d,m)\in \mathbb{Q}[m][d]$. After clearing denominators, let
$\widetilde{H}(d,m)\in \mathbb{Z}[d,m]$ be the corresponding polynomial.
The reduction of $\widetilde{H}(d,m)$ modulo $d$ is a nonzero constant
$c$ independent of $m$. Hence, if $\widetilde{H}(d,m)=0$ for some integer
$m$, then $d$ divides $c$. This is impossible for all sufficiently large
$d$. Therefore,
\[
\operatorname{Res}_x(P(x),g(x))\ne 0
\]
for all sufficiently large $d$.

On the other hand, both $P(x)$ and $g(x)$ vanish at $\alpha$. This
contradicts the nonvanishing of the resultant. Therefore, almost tight
spherical $2e$-designs in $S^{d-1}$ do not exist for sufficiently large
$d$.
\end{proof}

\subsection{Classification of almost tight spherical $2e$-designs for $e\leq 12$}

From the argument in the previous subsection, we can write
\[
\operatorname{Res}_x(P(x),g(x))=d^2H_e(d,m)
\]
for some $H_e(d,m)\in\mathbb{Q}[d,m]$. 
Then, there exists a nonzero rational number $c_e$,
depending only on $e$, such that
\[
H_e(d,m)=c_e+dK_e(d,m)
\]
for some $K_e(d,m)\in\mathbb{Q}[d,m]$. 

Let $L_e$ be the smallest positive integer such that
\[
L_eH_e(d,m)\in\mathbb{Z}[d,m].
\]
Then, 
we have
\[
L_eH_e(d,m)\equiv L_ec_e \pmod d.
\]
Put
\begin{equation}\label{eq:defC_e}
C_e=L_ec_e.
\end{equation}
Then $C_e$ is a nonzero integer depending only on $e$.

If $d\geq 3e+1$ and an almost tight spherical $2e$-design exists, then
the argument in the previous subsection gives an integer $m$ such that
$\operatorname{Res}_x(P(x),g(x))=0$. Hence $H_e(d,m)=0$, and therefore
$d\mid C_e$.
Thus, for $d\geq 3e+1$, the only possible dimensions not excluded by
the resultant argument are the positive divisors of $C_e$. 
We can further eliminate many candidate dimensions using the necessary
integrality condition 
\[
|X_1|=\frac{1}{2}\left(\binom{d+e-1}{e}+\binom{d+e-2}{e-1}+1 \right) \in \mathbb{Z}.
\]

For the remaining candidates for $d$, we use the irreducibility of $P_e(x)$ and $Q_e(x)$ to prove the nonexistence. 
\begin{lemma} \label{lem:irr}
Let $e>1$ and $d>2$. Suppose that an almost tight spherical $2e$-design
$X$ exists, and let $X=X_1\sqcup X_2$ be the simplex preimage decomposition. 
\begin{enumerate}
    \item If $P_e(x)$ is irreducible over $\mathbb{Q}$, then 
    $|X_1|-1$ is divisible by $e$. 
    \item If $Q_e(x)$ is irreducible over $\mathbb{Q}$, then
    $|X_2|$ is divisible by $e$. 
\end{enumerate}
In particular, if both $P_e(x)$ and $Q_e(x)$ are irreducible over
$\mathbb{Q}$, then no almost tight spherical $2e$-design exists.
\end{lemma}

\begin{proof}
Let $\alpha_1,\ldots,\alpha_e$ and $\beta_1,\ldots,\beta_e$ be the zeros
of $P_e(x)$ and $Q_e(x)$, respectively. Since
$P_e(x)-Q_e(x)=2$, the polynomials $P_e(x)$ and $Q_e(x)$ have no common
zero.

Fix $x\in X_1$. For $1\leq i\leq e$, put
\[
v_i=
|\{y\in X_1\setminus\{x\}\mid \langle x,y\rangle=\alpha_i\}|.
\]
Define
\[
F_i(z)
=
\frac{P_e(z)Q_e(z)}
{(z-\alpha_i)P_e'(\alpha_i)Q_e(\alpha_i)}.
\]
Then $\deg F_i=2e-1$ and
\[
F_i(\alpha_j)=\delta_{ij},
\qquad
F_i(\beta_j)=0.
\]
Since $X$ is a spherical $2e$-design, we have
\[
F_i(1)+v_i
=
\sum_{y\in X}F_i(\langle x,y\rangle)
=
|X|(F_i)_0.
\]
Let $K$ be the splitting field of $P_e(x)Q_e(x)$ over $\mathbb{Q}$.
Since $P_e(x)$ is irreducible over $\mathbb{Q}$, the action of
$\operatorname{Gal}(K/\mathbb{Q})$ on the zeros of $P_e(x)$ is
transitive. 

Let $\sigma \in \operatorname{Gal}(K/\mathbb{Q})$ such that $\sigma(\alpha_i)=\alpha_j$. Then $\sigma(F_i)=F_j$. 
Since the Gegenbauer polynomials have rational coefficients, we have
\[
\sigma((F_i)_0)=(F_j)_0.
\]
Applying $\sigma$ to
\[
v_i=|X|(F_i)_0-F_i(1),
\]
we obtain $\sigma(v_i)=v_j$.
 Since $v_i$ is an integer, it is fixed by
$\sigma$, and hence $v_i=v_j$. Thus all $v_i$ are equal, say $v_i=a$. Therefore,
\[
|X_1|-1=\sum_{i=1}^e v_i=ae,
\]
and hence
\[
|X_1|\equiv1\pmod e.
\]

Similarly, for $1\leq i\leq e$, put
\[
w_i=
|\{y\in X_2\mid \langle x,y\rangle=\beta_i\}|.
\]
Using
\[
G_i(z)
=
\frac{P_e(z)Q_e(z)}
{(z-\beta_i)Q_e'(\beta_i)P_e(\beta_i)},
\]
the irreducibility of $Q_e(x)$ implies that all $w_i$ are equal, say
$w_i=b$. Hence
\[
|X_2|=\sum_{i=1}^e w_i=be,
\]
and therefore
\[
|X_2|\equiv0\pmod e.
\]
Since $|X_1|=|X_2|$, this is a contradiction.
\end{proof}

\begin{theorem}\label{thm:class12}
    An almost tight spherical $2e$-design in $S^{d-1}$ does not exist for any $d>2$ and $3\leq e \leq 12$.  
\end{theorem}
\begin{proof}
Fix $e\in\{3,\ldots,12\}$ and compute the integer $C_e$ defined in
\eqref{eq:defC_e}. Suppose that an almost tight spherical $2e$-design
exists with $d\geq 3e+1$. By the argument in the previous subsection, there
exists an integer $m$ such that
\[
\operatorname{Res}_x(P(x),g(x))=0.
\]
Hence $H_e(d,m)=0$, and therefore $d\mid C_e$.
Thus it remains to consider the finite set
\[
S_e=
\left\{
d\in\mathbb{Z}
\;\middle|\;
d\geq3e+1,\ 
d\mid C_e,\ 
\frac{1}{2}\left(
\binom{d+e-1}{e}
+
\binom{d+e-2}{e-1}
+1
\right)\in\mathbb{Z}
\right\}.
\]
By a direct computation with Mathematica and Python (see the Data Availability Statement), for every $d\in S_e$, both
$P_e(x)$ and $Q_e(x)$ are irreducible over $\mathbb{Q}$. Hence Lemma~\ref{lem:irr} proves the nonexistence for all $d\geq3e+1$. 
The values of $C_e$ and $|S_e|$ are given in Table~\ref{tab:1}.

It remains to consider $3\leq d<3e+1$. Among the dimensions for which
$|X_1|$ is an integer, both $P_e(x)$ and $Q_e(x)$ are irreducible over
$\mathbb{Q}$ except for
$
(e,d)=(3,6), (9,4), (12,4)$.
Thus all the other cases are excluded by Lemma~\ref{lem:irr}.

For $(e,d)=(3,6)$, the polynomial $P_3(x)$ is irreducible over
$\mathbb{Q}$. Hence, by Lemma~\ref{lem:irr} (1), $3$ must divide
$|X_1|-1=38$, a contradiction. 
For $(e,d)=(9,4)$, the polynomial $Q_9(x)$ is irreducible over
$\mathbb{Q}$. Hence, by Lemma~\ref{lem:irr} (2), $9$ must divide
$|X_2|=193$, a contradiction. 
For $(e,d)=(12,4)$, the polynomial $P_{12}(x)$ is irreducible over
$\mathbb{Q}$. Hence, by Lemma~\ref{lem:irr} (1), $12$ must divide
$|X_1|-1=409$, a contradiction.

Therefore, no almost tight spherical $2e$-design in $S^{d-1}$ exists for
$d>2$ and $3\leq e\leq12$.
\end{proof}
\begin{table}[htbp]
\centering
\begin{tabular}{c|ccccc}
$e$ & 3 & 4 & 5 & 6 & 7 \\ \hline
$C_e$
& $2^{17}3^3$
& $2^{37}3^{10}$
& $2^{65}3^{12}5^7$
& $2^{92}3^{23}5^{16}$
& $2^{122}3^{35}5^{18}7^{11}$ \\
$|S_e|$
& $2$
& $18$
& $156$
& $308$
& $4101$
\end{tabular}

\medskip

\begin{tabular}{c|ccc}
$e$ & 8 & 9 & 10 \\ \hline
$C_e$
& $2^{185}3^{39}5^{20}7^{22}$
& $2^{247}3^{73}5^{22}7^{24}$
& $2^{298}3^{103}5^{41}7^{26}$ \\
$|S_e|$
& $48308$
& $85130$
& $88447$
\end{tabular}

\medskip

\begin{tabular}{c|cc}
$e$ & 11 & 12 \\ \hline
$C_e$
& $2^{348}3^{111}5^{57}7^{28}11^{19}$
& $2^{430}3^{140}5^{61}7^{30}11^{34}$ \\
$|S_e|$
& $1883836$
& $9485092$
\end{tabular}

\caption{The values of $C_e$ and $|S_e|$ for $3\leq e\leq12$.}
\label{tab:1}
\end{table}

\section{Spherical $2e$-designs with larger cardinalities}\label{sect:dual}

Our results on almost tight spherical designs rely on the classification of $2$-designs of size $d+2$ in $S^{d-1}$ found by Mimura \cite{M90} (with the uniqueness up to orthogonal transformations shown in \cite{S93}). $2$-designs of larger size cannot be classified in the same manner. Moreover, they need not be rigid as designs \cite{S93} (see the definition of rigidity in \cite{B87}). However, the approach of Mimura still provides useful structural information about designs. The following theorem is essentially a more elaborate restatement of \cite[Section 2]{M90}.

\begin{theorem}\label{thm:2-designs}
Let $X$ be a $2$-design in $S^{d-1}$ such that $|X|=d+k$, $k\geq 2$. Then there exists a $2$-design $Y$ in $S^{k-2}$ such that $|Y|=d+k$ and a one-to-one correspondence $x_i\leftrightarrow y_i$ between the elements of multisets $X=\{x_1,\ldots, x_{d+k}\}$ and $Y=\{y_1, \ldots, y_{d+k}\}$ such that $\langle x_i, x_j \rangle = -\frac 1 d -\frac {k-1} d \langle y_i, y_j \rangle$ for all pairs $i\neq j$.
\end{theorem}

\begin{proof}
Let $X=\{x_1,\ldots, x_{d+k}\}$ be a $2$-design in $S^{d-1}$. Using the design condition from Theorem~\ref{thm:equivalent_design}, we obtain
\[
\sum\limits_{i=1}^{d+k} \langle z, x_i\rangle = 0,\qquad \sum\limits_{i=1}^{d+k} \langle z, x_i\rangle^2 = \frac{d+k}{d}\| z \|^2
\]
for every $z\in \R^d$. Taking $z=z_1+z_2$ in the second identity, we also obtain
\[
\sum\limits_{i=1}^{d+k} \langle z_1, x_i\rangle \langle z_2,x_i\rangle = \frac{d+k}{d}\langle z_1, z_2 \rangle
\]
for arbitrary vectors $z_1$ and $z_2$.

Now we construct a matrix of size $d\times (d+k)$ whose columns are column vectors $\sqrt{\frac d {d+k}} x_i$. Then we add one more row to this matrix all of whose entries are $\sqrt{\frac 1 {d+k}}$. Using the three identities above with vectors $z$, $z_1$, $z_2$ from the standard basis, we can see that this $(d+1)\times (d+k)$ matrix has orthonormal rows.

We consider an arbitrary completion of the matrix to an orthogonal matrix of order $d+k$. Each column is extended by a vector in $\R^{k-1}$ with norm $\sqrt{\frac {k-1} {d+k}}$. By $y_i$ we denote such an extension of column $i$ multiplied by $\sqrt{\frac {d+k} {k-1}}$ so that all vectors $y_i$ are unit vectors in $\R^{k-1}$.

The matrix is orthogonal so for any distinct $i$ and $j$,
\[
\frac d {d+k} \langle x_i, x_j \rangle + \frac 1 {d+k} + \frac {k-1} {d+k} \langle y_i, y_j \rangle = 0,
\]
which confirms the connection between inner products in the statement of the theorem.

The set $Y=\{y_1,\ldots,y_{d+k}\}$ is a $2$-design because it satisfies the same three identities above so we can use Theorem~\ref{thm:equivalent_design}. These identities are automatically satisfied when vectors $z$, $z_1$, $z_2$ are from the standard basis due to the orthogonality of the matrix. For any other vectors, we can see that the identities hold using their expansions in the standard basis.
\end{proof}

\begin{remark}
Theorem~\ref{thm:2-designs} essentially establishes duality between $2$-designs of the same size $d+k$ in dimensions $d$ and $k-1$. The same kind of duality is very well known for tight frames (also called spherical designs of harmonic index $2$ when all vectors in a frame are unit). In particular, for every unit norm tight frame of size $d+k$ in $\R^d$, there exists a so-called \textit{Naimark dual} unit norm tight frame of size $d+k$ in $\R^k$ \cite{N43}. The duality in Theorem~5.1 can also be formulated using the complements of balanced frames introduced in~\cite[Section~6]{HMT20}.

\end{remark}

Theorem~\ref{thm:harmmap} provides the construction of a $2$-design based on a $2e$-design. As such, we can combine it with Theorem~\ref{thm:2-designs} to find certain restrictions on spherical $2e$-designs.

\begin{theorem}\label{thm:2e-dual}
Let $X=\{x_1, \ldots, x_N\}$ be a spherical $2e$-design in $S^{d-1}$ of size $N\geq M(d,e)+1$. Then there exists a spherical $2$-design $Y=\{y_1,\ldots,y_N\}$ in $S^{N-M(d,e)-1}$ and a one-to-one correspondence $x_i\leftrightarrow y_i$ between the multisets $X$ and $Y$ such that
\[
R_e(\langle x_i, x_j \rangle) =  -(N-M(d,e)) \langle y_i, y_j\rangle \quad\forall i\neq j.
\]
\end{theorem}

\begin{proof}
Let $H=h_1+\ldots+h_e=M(d,e)-1$ as above. We first consider $\Phi(X)$ in $S^{H-1}$ that is a $2$-design by Theorem~\ref{thm:harmmap}. Pairwise inner products between the points of $\Phi(X)$ are precisely $\frac {R_e(\langle x_i, x_j\rangle)-1} H$. Now we apply Theorem~\ref{thm:2-designs} to $\Phi(X)$ and construct a $2$-design $Y=\{y_1,\ldots,y_N\}$ in $S^{N-H-2}$ whose inner products satisfy
\[
\frac {R_e(\langle x_i, x_j\rangle)-1} H = -\frac 1 H - \frac {N-H-1} {H} \langle y_i, y_j \rangle \quad i\neq j,
\]
from which the required correspondence follows immediately.
\end{proof}

\begin{remark}\label{rem:tight}
We note that the proof of the lower bound for the size of a $2e$-design by Delsarte, Goethals, and Seidel was based precisely on polynomial $R_e$ and $R_e(1)=M(d,e)$ \cite{DGS77}. Reproducing the proofs of Theorem~\ref{thm:2-designs} for $k=1$ and Theorem~\ref{thm:2e-dual} for $N=M(d,e)$ we get that $\langle x_i,x_j \rangle$ are zeros of the polynomial $R_e$ which is one of the main properties used in a series of papers by Bannai and coauthors to classify tight designs \cite{BD79, BD80}.
\end{remark}

The case of almost tight designs is covered in this paper in a somewhat similar manner because the dual design $Y$ is one-dimensional and all its inner products are $\pm 1$ which means that all pairwise inner products $\langle x_i, x_j\rangle$ can be found explicitly as the zeros of the polynomials $R_e\pm 1$.

There are no similar constraints when $N\geq M(d,e)+2$. We cannot find inner products of a $2e$-design $X$ of size $N$ explicitly. However, we can still use that $\langle y_i, y_j\rangle \in [-1,1]$ so Theorem~\ref{thm:2e-dual} allows us to obtain constraints on pairwise inner products of $X$.

\begin{theorem}\label{thm:scalar_prod}
Let $N\geq M(d,e)$. Let $I:=I(N,d,e)$ be the set $\{t\in[-1,1]: |R_e(t)|\leq N-M(d,e)\}$. If $X$ is a $2e$-design in $S^{d-1}$ of size $N$, then $A(X)\subset I$.
\end{theorem}
\begin{proof}
This follows immediately from Theorem~\ref{thm:2e-dual} for $N\geq M(d,e)+1$ and Remark~\ref{rem:tight} for $N=M(d,e)$.
\end{proof}

We note that this theorem does not provide any new information about designs if $N\geq 2M(d,e)$ because $|R_e(t)|\leq M(d,e)$ for all $t\in [-1,1]$. For all smaller values of $N$ it is a new nontrivial constraint on inner products. It will be used in the next section to set up a semidefinite program on inner products of designs.

\section{Semidefinite programs for designs of nearly tight size}\label{sect:3-point}

In this section, we develop a 3-point semidefinite program that allows us to obtain polynomial certificates of nonexistence for $2e$-designs of certain sizes $<2M(d,e)$, that is, for those where Theorem~\ref{thm:scalar_prod} provides a nontrivial set of restrictions $I(N,d,e)$ on inner products in a design. The program is based on a 3-point semidefinite program by Bachoc and Vallentin \cite{BV08} with several relaxations and modifications: there are no positive semidefinite constraints on harmonic degrees $1,\ldots,2e$ because of the design conditions; all inner products are from $I$ (Theorem~\ref{thm:scalar_prod}); the Gram matrix of all triples of points in the dual configuration from Theorem~\ref{thm:2e-dual} is positive semidefinite. In the case $N=M(d,e)+2$ the dual configuration lies on a circle, so we use an even stronger constraint that every $3\times3$ Gram determinant vanishes.

\subsection{Bachoc-Vallentin semidefinite program}\label{subsect:BV}

Here we give a short overview of the 3-point SDP by Bachoc and Vallentin \cite{BV08}.

For points $x,y,z$ in the sphere $S^{d-1}$ we denote $u=\langle x,z \rangle$, $v=\langle y,z\rangle$, $t=\langle x,y\rangle$. By $P^{(d)}_{i}$ we denote the Gegenbauer polynomial of degree $i$ normalized at $1$ for the dimension $d$, that is, 
$P^{(d)}_i(u)=G^{(d)}_{i}(u)/G^{(d)}_{i}(1).$ Define, following~\cite[Section~3]{BV08},
\begin{equation}\label{eq:Q}
Q^{d-1}_k(u,v,t)=\bigl((1-u^2)(1-v^2)\bigr)^{k/2}\,P^{d-1}_k\!\left(\frac{t-uv}{\sqrt{(1-u^2)(1-v^2)}}\right),
\end{equation}
a polynomial in $(u,v,t)$ of degree $k$ in each of the variables.

Suppose $\xi$ is a fixed point in $S^{d-1}$. Let $Y^d_k\in\R[u,v,t]^{(D-k+1)\times(D-k+1)}$ be a matrix-valued polynomial function on $[-1,1]^3$ (a matrix of
polynomials) defined by
\begin{equation}\label{eq:Y}
\bigl(Y^d_k\bigr)_{ij}(u,v,t)=\lambda_{ij}\,P^{d+2k}_i(u)\,P^{d+2k}_j(v)\,Q^{d-1}_k(u,v,t),\qquad 0\le i,j\le D-k,
\end{equation}
where $\lambda_{ij}$ are suitably chosen positive constants and $D$ is the truncation parameter. Let $Z^d_k(x,y)$ be the value of $Y^d_k$ at $(\xi\cdot x,\xi\cdot y,x\cdot y)$. Then for every finite
$C\subset S^{d-1}$, $\sum_{(c,c')\in C^2}Z^d_k(c,c')\succeq 0$ \cite[Theorem~3.1]{BV08}. The proof of Bachoc and Vallentin is based on irreducible representations of the stabilizer of $\xi$, as we will explain below, but a somewhat simpler explanation just uses positive semidefiniteness of Gegenbauer polynomials on the equator \cite{S42} and applies it to equatorial projections of the points of $C$ onto $\xi^{\perp}$.

Finally, \cite[Corollary~3.5]{BV08} introduces the symmetrization
\[
S^d_k(u,v,t)=\frac16\sum_{\sigma\in\mathfrak S_3}Y^d_k\bigl(\sigma(u,v,t)\bigr)\in\R[u,v,t]^{(D-k+1)\times(D-k+1)},
\]
where $\mathfrak S_3$ acts by permuting the variables $(u,v,t)$ in every entry so that $S^d_k(u,v,t)$ is a symmetric matrix for every $(u,v,t)$. The key theorem for the SDP is the following.

\begin{theorem}\cite[Corollary~3.5]{BV08}\label{thm:BV}
For every $(u,v,t)$ the matrix $S^d_k(u,v,t)$ is symmetric, its entries are symmetric polynomials in $(u,v,t)$, and
for every finite $C\subset S^{d-1}$ and every $k\ge 0$,
\[
\sum_{(c,c',c'')\in C^3}S^d_k\bigl(c\cdot c',\,c\cdot c'',\,c'\cdot c''\bigr)\ \succeq\ 0 .
\]
\end{theorem}

The proof in~\cite{BV08} is short: by the positive semidefiniteness of $Z_k^d$ applied with the pole $\xi=c''\in C$ one has
$\sum_{(c,c')\in C^2}Y^d_k(c''\cdot c,\,c''\cdot c',\,c\cdot c')\succeq 0$ for each $c''$; summing over $c''\in C$ gives the statement for the unsymmetrized $Y^d_k$, and the aggregate over $C^3$ is invariant under $\mathfrak S_3$, so it coincides with the aggregate of $S^d_k$.

\subsection{Design relaxations of the 3-point SDP}\label{subsect:relax}

Here we show which conditions of the SDP from Subsection~\ref{subsect:BV} can be relaxed in the case of $2e$-designs. First, we show how the rows of $Y^d_k$ are built from spherical harmonics. The arguments below are directly taken from \cite[Section 3]{BV08}.

As above, $\xi\in S^{d-1}$ is a fixed point. Its stabilizer is $H\cong O(d-1)$. Decompose $\operatorname{Pol}_{\le D}(S^{d-1})$ into $H$-isotypic components $\mathcal I_k\cong(D-k+1)\operatorname{Harm}_k(\R^{d-1})$, $0\le k\le D$. Then each $\mathcal I_k$ corresponds to the matrix-valued zonal function $Z^d_k(x,y)=E(x)E(y)^{\mathsf T}$ we defined above, where the rows of $E(x)$ are orthonormal bases of the $D-k+1$ irreducible constituents of $\mathcal I_k$. The positive semidefiniteness of $Z^d_k$ follows from this right away.

In the proof of \cite[Theorem~3.2]{BV08} the rows of $E(x)$ are made explicit: writing $x=u\xi+\sqrt{1-u^2}\,\zeta$ with
$\zeta\in S^{d-2}$, the $i$-th row of $E(x)$ consists of the functions
\begin{equation}\label{eq:harmdeg}
e^{(\xi)}_{k,i,m}(x)=\lambda_i\,P^{d+2k}_i(u)\,(1-u^2)^{k/2}f_m(\zeta),\qquad m=1,\dots,h^{d-1}_k,
\end{equation}
where $h^{d-1}_k$ is the dimension of $\operatorname{Harm}_k(\mathbb{R}^{d-1})$, $(f_m)$ is an orthonormal basis of $\operatorname{Harm}_k(\mathbb{R}^{d-1})$, $\lambda_i>0$ is a normalizing constant, and the superscript records the pole $\xi$ (we write $e^{(z)}_{k,i,m}$ when the pole is a point $z\in C$). These functions form an orthonormal basis of the subspace $H^{d-1}_{k,k+i}$ of $H^d_{k+i}\cong\operatorname{Harm}_{k+i}(\mathbb{R}^d)$, the constituent of $\operatorname{Harm}_{k+i}(\mathbb{R}^d)$ isomorphic to $\operatorname{Harm}_k(\mathbb{R}^{d-1})$ under $H$ (see also \cite[Ch.~9.8]{AAR99}). In other words, row $i$ of the matrix $E$ defining $Z^d_k$ is built from spherical harmonics of degree exactly $k+i$.

Let $C$ be a spherical $2e$-design with $N$ points and let $D$ be the degree parameter. By Theorem~\ref{thm:equivalent_design}, we can describe two-point relaxations for $2e$-designs:
\[
\sum_{x,y\in C}P^d_k(\langle x,y\rangle)=
\begin{cases}N^2,&k=0,\\ 0,&1\le k\le 2e,\\ \ge 0,&k>2e.\end{cases}
\]
In a similar manner we can now describe design relaxations for three-point constraints.
\begin{lemma}\label{lem:threepoint}
Let $M_k:=\sum_{(x,y,z)\in C^3}S^d_k(\langle x,z\rangle,\langle y,z\rangle, \langle x,y\rangle)$, $0\le k\le D$. Then row and column $i$ of $M_k$ vanish whenever $1\le k+i\le 2e$, and the principal submatrix of $M_k$ indexed by the \emph{active set}
\[
A_k:=\{\,i\in\{0,\dots,D-k\}:\ k+i=0\ \text{or}\ k+i>2e\,\}
\]
is positive semidefinite.
\end{lemma}
\begin{proof}
By \cite[Theorems~3.1 and~3.2]{BV08} (as described in Section~\ref{subsect:BV}), with the pole $\xi=z$,
\[
\sum_{(x,y,z)\in C^3}(Y^d_k)_{ij}=c_k\sum_{z\in C}\sum_{m}\Bigl(\sum_{x\in C}e^{(z)}_{k,i,m}(x)\Bigr)
\Bigl(\sum_{y\in C}e^{(z)}_{k,j,m}(y)\Bigr),\qquad c_k>0,
\]
where $e^{(z)}_{k,i,m}$ is the function~\eqref{eq:harmdeg} for the pole $z$, a spherical harmonic of degree $k+i$ on
$S^{d-1}$. Since $C$ is a $2e$-design, $\sum_{x\in C}e^{(z)}_{k,i,m}(x)=0$ for every $z\in C$ whenever $1\le k+i\le 2e$;
hence row and column $i$ of $\sum_{C^3}Y^d_k$ vanish for those $i$, and the remaining principal submatrix is a sum of Gram matrices, hence positive semidefinite. The symmetrized aggregate $M_k$ equals the unsymmetrized one because the aggregate over $C^3$ is invariant under permutations of the three points.
\end{proof}

Bachoc and Vallentin remark that the polynomials $P^{d+2k}_i(u)$ may be replaced by the monomials $u^i$ which may be more convenient for computations, but the harmonic-degree property holds only after an invertible recombination of the rows, which is why we keep the basis as defined above. For the sake of convenience, we drop the positive constants $\lambda_{ij}$: replacing $Y^d_k$ by $\Lambda^{-1}Y^d_k\Lambda^{-1}$
with $\Lambda=\operatorname{diag}(\lambda_{ii}^{1/2})$ is a positive diagonal conjugation, which preserves positive semidefiniteness of every aggregate and of every principal submatrix. From this point onward, $Y^d_k$ and $S^d_k$ denote these rescaled matrix-valued polynomials, with entries $P^{d+2k}_i(u)P^{d+2k}_j(v)Q^{d-1}_k(u,v,t)$ and their
symmetrizations; for a symmetric matrix $F_k\in\mathcal{S}^{D-k+1}(\mathbb{R})$ the expression $\langle F_k,S^d_k(u,v,t)\rangle=\operatorname{tr}\bigl(F_kS^d_k(u,v,t)\bigr)$ is a symmetric polynomial in $(u,v,t)$. The conclusions of Lemma~\ref{lem:threepoint} remain valid for these rescaled kernels.

\begin{corollary}\label{cor:relax}
Let $F_k\in\mathcal{S}^{D-k+1}(\R)$, $0\le k\le D$, be symmetric matrices such that the principal submatrix $F_k[A_k,A_k]$ is positive semidefinite; the entries of $F_k$ outside $A_k\times A_k$ are arbitrary. Then
\begin{equation}\label{eq:Fpos}
\sum_{(x,y,z)\in C^3}F\bigl(\langle x,z\rangle,\langle y,z\rangle,\langle x,y\rangle\bigr)\ \ge\ 0,\qquad
F(u,v,t):=\sum_{k=0}^{D}\langle F_k,S^d_k(u,v,t)\rangle.
\end{equation}
\end{corollary}
\begin{proof}
$\sum_{C^3}F=\sum_k\langle F_k,M_k\rangle=\sum_k\langle F_k[A_k,A_k],M_k[A_k,A_k]\rangle\ge 0$ by Lemma~\ref{lem:threepoint}.
\end{proof}

\begin{remark}
When $D=2e$ the only active entry is $(F_0)_{00}$ (harmonic degree $0$); every other entry of every $F_k$ is free. Apart from the condition $(F_0)_{00}\ge 0$, Corollary~\ref{cor:relax} then uses no positive semidefiniteness: it is the system of linear identities expressing that the aggregates of all entries of harmonic degree $1,\dots,2e$ vanish, and the program of Subsection~\ref{subsect:SDP} is a semi-infinite linear program in the three-point function $F$. This is the case for $226$ of the $265$ certificates of Subsection~\ref{subsect:certs}, namely for all those with $(e,D)=(3,6)$, $(4,8)$ and $(5,10)$. In each of them $(F_0)_{00}$ is zero or negligible: it is exactly $0$ in $85$ certificates and at most $2\cdot 10^{-6}$ in the others. The remaining $39$ certificates, among them all certificates for $e=2$, have $D>2e$, and their positive semidefinite blocks $F_k[A_k,A_k]$ are far from zero.
\end{remark}

\subsection{Semidefinite program}\label{subsect:SDP}

Here we first determine the admissible domain for triples of inner products of three points. Then we formulate the semidefinite program that can produce certificates of nonexistence for designs of a given size.

From now on $C$ is a $2e$-design with $N=M+s$ points, $s\ge 1$, where $M=M(d,e)$, the size of a tight design. 

\begin{proposition}\label{prop:domain}
Let $C$ be a spherical $2e$-design with $N=M(d,e)+s$ points, $s\ge 1$. Let $P(t)=-R_e(t)/s$. Then
\begin{enumerate}
\item[\textup{(I)}] every inner product $u=\langle x,x'\rangle$, $x\neq x'$, lies in
$I:=\{t\in[-1,1]:|R_e(t)|\le s\}$, a union of at most $e$ closed intervals contained in $[-1,1]$;
\item[\textup{(g)}] $g(u,v,t):=1+2uvt-u^2-v^2-t^2\ge 0$ for every triple of distinct points;
\item[\textup{(h)}] $h(u,v,t):=1+2P(u)P(v)P(t)-P(u)^2-P(v)^2-P(t)^2\ge 0$ for every triple of distinct points, with
\begin{itemize}
\item $h(u,v,t)=0$ if $s=2$;
\item if $s=1$, $P(u)\in\{\pm1\}$ for every pair and $P(u)P(v)P(t)=1$ for every triple.
\end{itemize}
\end{enumerate}
\end{proposition}

\begin{proof}
Property (I) is just Theorem~\ref{thm:scalar_prod}. Property (g) follows from the positive semidefiniteness of a $3\times 3$ Gram submatrix. Property (h) follows from the nonnegativity of the determinant of a $3\times 3$ Gram submatrix in the dual construction of Theorem~\ref{thm:2e-dual}. When $s=2$, the dual construction lives in dimension 2 so any $3\times 3$ principal Gram submatrix has determinant 0. When $s=1$, the dual construction is in $\{-1,1\}$.
\end{proof}

We write $\mathcal D$ for the set of triples $(u,v,t)\in I^3$ satisfying (g) and (h) (with $h=0$ when $s=2$, and the finite set described above when $s=1$), so that every triple of distinct points of $C$ has its inner products in $\mathcal D$. The constraints (I) and (h) are those that an arbitrary spherical code does not satisfy; they encode the design property of $C$.

We are ready to formulate the semidefinite program.

\begin{definition}[SDP for designs]\label{def:des_SDP}
Let $d,e,s,D$ be given, $N=M(d,e)+s$, $P(t)=-R_e(t)/s$, and $I,\mathcal D$ as in Proposition~\ref{prop:domain}. The variables are
\begin{itemize}
\item real numbers $a_0,a_1,\dots,a_{2D}$ with
\begin{equation}\label{eq:signs}
a_0=0,\qquad a_1,\dots,a_{2e}\ \text{free},\qquad a_k\ge 0\ \ (2e<k\le 2D),
\end{equation}
defining $f(t)=\sum_{k=0}^{2D}a_kP^d_k(t)$;
\item symmetric matrices $F_k\in\mathcal{S}^{D-k+1}(\mathbb{R})$, $0\le k\le D$, with $F_k[A_k,A_k]\succeq 0$ (Lemma~\ref{lem:threepoint}), defining $F(u,v,t)=\sum_k \langle F_k,S^d_k(u,v,t)\rangle$.
\end{itemize}
The constraints are
\begin{align}
f(u)+3F(u,u,1)&\le -1\qquad\text{for all }u\in I,\label{eq:B1}\\
F(u,v,t)&\le 0\qquad\ \ \text{for all }(u,v,t)\in\mathcal D.\label{eq:B2}
\end{align}
The objective is to minimize
\begin{equation}\label{eq:bound}
B:=1+f(1)+F(1,1,1),\qquad f(1)=\sum_{k}a_k,\quad F(1,1,1)=\sum_{i,j}(F_0)_{ij}
\end{equation}
(note $Q^{d-1}_k(1,1,1)=0$ for $k\ge 1$, so only $F_0$ contributes to $F(1,1,1)$).
\end{definition}

\begin{theorem}\label{thm:main}
If $C\subset S^{d-1}$ is a spherical $2e$-design with $N=M(d,e)+s$ points and $(a,F)$ is feasible for the semidefinite program described above,
then $N\le B$. Consequently, if the optimal value of the program is $<N$, or more simply if the feasibility problem $\{$\eqref{eq:signs}, $F_k[A_k,A_k]\succeq 0$, \eqref{eq:B1}, \eqref{eq:B2}, $B\le B_0\}$, where $B_0<N$ is a fixed number, has a solution, then no such design exists.
\end{theorem}

\begin{proof}
Split the sum in Corollary~\ref{cor:relax} according to coincidences among $x,y,z$: the $N$ terms with $x=y=z$ give $NF(1,1,1)$; the $3N(N-1)$ terms with exactly two equal points give $3\sum_{x\ne y}F(u_{xy},u_{xy},1)$ (by symmetry of $F$ all three coincidence patterns contribute equally); the remaining terms are over triples of
distinct points. Hence
\begin{equation}\label{eq:split}
0\le NF(1,1,1)+3\sum_{x\neq y}F(u_{xy},u_{xy},1)+\sum_{x,y,z\ \text{distinct}}F(u_{xz},u_{yz},u_{xy}).
\end{equation}
By Proposition~\ref{prop:domain} every triple of distinct points has inner products in $\mathcal D$, so the last sum is
$\le 0$ by \eqref{eq:B2}. Every $u_{xy}$ lies in $I$, so by \eqref{eq:B1}
$3F(u_{xy},u_{xy},1)\le -1-f(u_{xy})$, whence
\[
0\le NF(1,1,1)-N(N-1)-\sum_{x\neq y}f(u_{xy}).
\]
Finally, by the two-point conditions for designs and \eqref{eq:signs},
\[
\sum_{x\neq y}f(u_{xy})=\sum_{k}a_k\Bigl(\sum_{x,y\in C}P^d_k(u_{xy})-N\Bigr)=
a_0(N^2-N)+\sum_{k>2e}a_k\sum_{x,y}P^d_k(u_{xy})-N\sum_{k\ge1}a_k\ \ge\ -Nf(1),
\]
because $a_0=0$ and $a_k\ge0$, $\sum_{x,y}P^d_k\ge0$ for $k>2e$. Combining, $0\le NF(1,1,1)-N(N-1)+Nf(1)$, i.e.\
$N-1\le F(1,1,1)+f(1)$, which is $N\le B$.
\end{proof}

\subsection{Numerical solution and rounding}\label{subsect:numer}

The constraints \eqref{eq:B1}--\eqref{eq:B2} are polynomial inequalities on semialgebraic sets. We replace them, exactly as in \cite[Section~5]{BV08}, by stronger constraints through weighted sums-of-squares (SOS) identities in the spirit of the Positivstellens\"atze of Schm\"udgen~\cite{S91} and Putinar~\cite{P93} and the SOS relaxations of
Lasserre~\cite{L01} and Parrilo~\cite{P03}. A polynomial $\sigma$ is SOS if $\sigma=m^{\mathsf T}Gm$ for a vector $m$ of monomials and a positive semidefinite Gram matrix $G$. Write
\[
w(t):=1-P(t)^2=\bigl(1-P(t)\bigr)\bigl(1+P(t)\bigr),
\]
where $P(t)=-R_e(t)/s$ as above. Then $w(t)$ is a polynomial of degree $2e$ which is nonnegative exactly on $I$ for $t\in [-1,1]$. We use the following SOS form of the program.

\begin{definition}[SOS program for designs]\label{def:sos}
Find $a$, $(F_k)$ as in Definition~\ref{def:des_SDP}, SOS polynomials $\sigma_0,\sigma_1$ in $t$, SOS polynomials
$\tau_0,\tau_u,\tau_v,\tau_t,\tau_g$ in $(u,v,t)$, and a multiplier $\lambda(u,v,t)$, such that
\begin{align}
-1-f(t)-3F(t,t,1)&=\sigma_0(t)+w(t)\,\sigma_1(t),\label{eq:B1sos}\\
-F(u,v,t)&=\tau_0+w(u)\tau_u+w(v)\tau_v+w(t)\tau_t+g\,\tau_g+h\,\lambda ,\label{eq:B2sos}
\end{align}
as identities between polynomials, where
\begin{itemize}
\item for $s=2$: $\lambda$ is an arbitrary polynomial ($h=0$ on $\mathcal D$);
\item for $s\ge3$: $\lambda$ is SOS ($h\ge0$ on $\mathcal D$);
\item for $s=1$: the multipliers of $w(u),w(v),w(t)$ are arbitrary polynomials and $\lambda$ is an arbitrary polynomial multiplying $P(u)P(v)P(t)-1$ instead of $h$.
\end{itemize}
\end{definition}

Clearly \eqref{eq:B1sos} implies \eqref{eq:B1} and \eqref{eq:B2sos} implies \eqref{eq:B2}, so every solution of Definition~\ref{def:sos} yields a feasible solution of Definition~\ref{def:des_SDP} and Theorem~\ref{thm:main} applies. In all our computations both sides of \eqref{eq:B1sos} and \eqref{eq:B2sos} have degree at most $2D$: the Gram matrices of $\sigma_0$ and $\tau_0$ are indexed by all monomials of degree at most $D$, and every other multiplier has the largest degree for which its term has degree at most $2D$.

The SOS form is a semidefinite program with free variables ($a_1,\dots,a_{2e}$, the inactive entries of the $F_k$,
the coefficients of $\lambda$ for $s=2$), nonnegative variables ($a_k$, $k>2e$) and positive semidefinite blocks
($F_k[A_k,A_k]$ and the Gram matrices for SOS polynomials). We solve it as a feasibility problem with the additional constraint $B\le B_0$ for a chosen target $B_0<N$. Most certificates were obtained with the double-precision interior-point solver Clarabel~\cite{GC26}. When double precision was not accurate enough for the rounding below, in particular for all certificates with $e=5$, we used the multiprecision solver SDPA-GMP~\cite{N10}.

Following the general strategy of Peyrl and Parrilo~\cite{PP08} (see also~\cite{BV08, DLM21}), a numerical solution is converted into an exact rational one.
\begin{enumerate}
\item Choose a common denominator $\mathrm{den}$ ($10^{16}$ for most certificates). Round $a_k$, the free entries of the $F_k$, and the
coefficients of $\lambda$ for $s=2$ to $\mathrm{den}^{-1}\mathbb{Z}$; set $a_0=0$ and replace negative $a_k$, $k>2e$, by $0$.
\item For every positive semidefinite block $X$ compute an integer matrix $L$ with
$X\approx LL^{\mathsf T}/\mathrm{den}^2$ (a rounded factor of $X$) and define the exact block as
$LL^{\mathsf T}/\mathrm{den}^2$, which is positive semidefinite by construction.
\item The Gram matrices $G_{\sigma_0}$ and $G_{\tau_0}$ of the SOS terms without a multiplier are then corrected so that \eqref{eq:B1sos}--\eqref{eq:B2sos} hold exactly. To leave room for the correction, the numerical problem is solved with the
constraint $G\succeq\varepsilon W$ for these two matrices, where $\varepsilon>0$ is a small rational number (between $10^{-9}$ and $0.03$) and $W=\operatorname{diag}(v_i^2)$, where $v_i$ are fixed positive rational numbers. The exact matrix is
$G=LL^{\mathsf T}/\mathrm{den}^2+\varepsilon W+\Delta$, where the symmetric rational matrix $\Delta$ absorbs the residual of the identity
after all other terms have been rounded: the coefficient of each monomial $\mu$ of the residual is distributed over the entries $(i, j)$ with $m_i m_j = \mu$ in proportion to $v_i v_j$, where $m_i$ is the $i$-th monomial of the basis.
\item The certificate is valid if $G_{\sigma_0}\succeq 0$, $G_{\tau_0}\succeq 0$ and $B<N$. For $220$ of the $265$ certificates positive semidefiniteness follows from $\varepsilon W+\Delta\succeq0$, which is verified by the weighted Gershgorin
criterion: for every row $i$, $\varepsilon W_{ii}+\Delta_{ii}\ge\sum_{j\ne i}|\Delta_{ij}|\sqrt{W_{ii}/W_{jj}}$. For the other $45$ we verify $G\succeq 0$ directly: $W^{-1/2}GW^{-1/2}$ is written as $RR^{\mathsf T}+E$, where $R$ is a rounded floating-point Cholesky factor and the remainder $E$, computed exactly, is diagonally dominant with nonnegative diagonal.
\end{enumerate}
All checks in (4) are exact rational computations and do not depend on how the certificate was found. The certificates and a short program that verifies them are publicly available; see the Data Availability Statement.

\begin{table}[ht!]\centering
\caption{New lower bounds for the size of spherical $2e$-designs in $S^{d-1}$: no such design has $M(d,e)+s'$ points for any $1\le s'\le s$, so every design that is not tight has at least $M(d,e)+s+1$ points. The size $M(d,e)+1$ is excluded by Theorems~\ref{thm:4-design} and~\ref{thm:class12}, and every other size by a certificate of Table~\ref{tab:certificates} or, in the three cells with $e=3$ and $d\in\{4,5,10\}$, by the linear programming bound (Table~\ref{tab:lp}).}\label{tab:bounds}
\begin{tabular}{@{}cccc@{}}
\toprule
$e$ & $d$ & $M(d,e)$ & $M(d,e)+s+1$ \\
\midrule
2 & 3 & 9 & 12 \\
2 & 4 & 14 & 17 \\
2 & 5 & 20 & 23 \\
\midrule
3 & 3 & 16 & 19 \\
3 & 4 & 30 & 34 \\
3 & 5 & 50 & 55 \\
3 & 10 & 275 & 285 \\
3 & 11 & 352 & 363 \\
3 & 12 & 442 & 454 \\
3 & 13 & 546 & 558 \\
3 & 14 & 665 & 677 \\
3 & 15 & 800 & 811 \\
3 & 16 & 952 & 960 \\
\midrule
4 & 3 & 25 & 28 \\
4 & 18 & 7125 & 7149 \\
4 & 19 & 8645 & 8669 \\
4 & 20 & 10395 & 10418 \\
4 & 21 & 12397 & 12421 \\
4 & 22 & 14674 & 14696 \\
4 & 23 & 17250 & 17271 \\
4 & 24 & 20150 & 20170 \\
4 & 25 & 23400 & 23418 \\
4 & 26 & 27027 & 27040 \\
\midrule
5 & 27 & 197316 & 197319 \\
5 & 28 & 232841 & 232872 \\
5 & 29 & 273296 & 273299 \\
\bottomrule
\end{tabular}
\end{table}

\subsection{Exact certificates of nonexistence}\label{subsect:certs}

For a polynomial $f=\sum_k a_kP^d_k$ with $a_0>0$, $a_k\le 0$ for $k>2e$, and $f(t)\ge 0$ on $[-1,1]$, every $2e$-design satisfies $N\ge f(1)/a_0$~\cite{DGS77}. Optimizing over $f$ gives the linear programming (LP) bound. Yudin~\cite{Y97} and Boyvalenkov and
Nikova~\cite{BN98} showed that it improves (\ref{eq:tight}) in some dimensions. The test functions of~\cite{BN98} improve (\ref{eq:tight}) for $2e$-designs and show nonexistence of tight designs precisely when $d\le e^2+1$ for small $e$. Nikova and Nikov~\cite{NN03} conjectured that for $d\ge e^2+2$ bound (\ref{eq:tight}) is the best bound obtainable by the linear programming method, that is, no admissible polynomial $f$ gives $f(1)/a_0>M$. This conjecture is open but heavily supported by numerical calculations. Our own LP computations (polynomials of degree up to $40$) found no improvement for $d\ge e^2+2$, consistently with the conjecture.

The three-point program described above excludes sizes $N=M+s$, $s\ge 2$, in the range in which no LP improvement on the Fisher-type bound is known. Below we provide a table of new lower bounds for sizes of $2e$-designs based on certificates of nonexistence found by the SOS program in Definition~\ref{def:sos}.  

All certificates are defined by a triple of parameters $(e,d,s)$. In all cases, the existence of an exact certificate implies nonexistence of a spherical $2e$-design in $S^{d-1}$ of size $M(d,e)+s$. We note that a separate certificate is required for each $s$ so, for example, a lower bound of 363 for 6-designs in dimension 11 required nine separate SDP certificates, one for each $s\in\{2,\ldots,10\}$. 

Our results can be split into two major groups: improvements to the lower bounds on cardinality of a spherical $2e$-design in the range where linear programming bounds are effective, that is, for $3\leq d\leq e^2+1$ and improvements on the size of a spherical $2e$-design beyond this range, that is, for $d\geq e^2+2$. In three cells of Table~\ref{tab:bounds}, namely $(e,d)=(3,4)$, $(3,5)$ and $(3,10)$, the smallest sizes are excluded already by the LP bound stated above: it gives $N\ge 33$, $N\ge 54$ and $N\ge 283$, respectively. Explicit polynomials are given in Table~\ref{tab:lp} of the appendix. The certificates of Table~\ref{tab:certificates} take over from there.

We use $D\ge 2e$ so that all harmonic degrees $1,\ldots,2e$ occur in the three-point kernels. With the computational resources of a personal laptop we can handle $D$ up to $12$, so the results are limited to small $e$. Most of the certificates are for $e=3$ and $e=4$, where $D=6$ or $D=8$ is sufficient.

At a fixed degree the method reaches only a limited range of dimensions. With $D=2e$ we found certificates for $11\le d\le 16$ when $e=3$, for $18\le d\le 26$ when $e=4$, and for $27\le d\le 29$ when $e=5$, that is, in a window that starts where the linear programming bound stops improving \eqref{eq:tight}. In the next dimensions ($d=17$ for $e=3$, $d=27$ for $e=4$, $d=30,31,32$ for $e=5$) we obtained no certificate: the solver reports infeasibility, does not converge, or returns a solution that cannot be rounded. To separate the limits of the program from those of the SOS relaxation and of the solver, we also solved relaxations of the program of Definition~\ref{def:des_SDP} in which \eqref{eq:B1} and \eqref{eq:B2} are imposed only on a finite sample of $I$ and of $\mathcal D$. Such a relaxation involves no SOS multipliers, and if its optimal value is at least $N$, then no certificate with the given $D$ exists. These computations are numerical and not rigorous, but they give a consistent picture. For $e=3$ and $s=2$ they indicate that no certificate with $D\le 8$ exists for $d=17$ and $d=20$, while $D=10$ may be sufficient for $d=17,20,25$ but not for $d=40$. These numerical results suggest that increasing the degree may extend the range of dimensions. We do not know whether, for fixed $e$ and $s$, every dimension is reached by some degree.

The same test shows why Table~\ref{tab:bounds} stops where it does in small dimensions. The sampled relaxation has optimal value greater than $N$ for $13$-point $4$-designs in $S^2$ for every $D\le 16$ (the value decreases from $15.01$ at $D=6$ to $13.55$ at $D=16$), for $17$-point $4$-designs in $S^3$ for every $D\le 12$, and for $34$-point $6$-designs in $S^3$ for every $D\le 12$. In particular, we cannot decide in this way whether $4$-designs with $13$ or $15$ points exist in $S^2$. According to the numerical results of Hardin and Sloane~\cite{HS92} they do not. Interestingly, it seems to be difficult to provide certificates of nonexistence for 4-designs in higher dimensions as they tend to require larger degrees. We do not have an explanation for this phenomenon.

Most notably, we found the SDP certificate of nonexistence for 4-designs of size 11 in $S^2$. Combined with Theorem~\ref{thm:4-design}, the nonexistence of a tight design, and the fact that the vertices of a regular icosahedron form a 5-design, this means that the smallest size of a 4-design in $S^2$ is 12. This and other new lower bounds on sizes of 4-designs are consistent with the numerical results of Hardin and Sloane \cite{HS92}. In Section~\ref{sect:constraint} we study constraints on inner products, particularly, for 12-point 4-designs in $S^2$ and find those with the maximal possible inner products in Theorem~\ref{thm:5-design_mod}.

Table~\ref{tab:certificates} in Appendix~\ref{app:certificates} lists the parameters and bounds for all certificates obtained.

\section{Constraints on distances in spherical designs}\label{sect:constraint}

In this section we study constraints on distances in spherical designs arising from Theorems~\ref{thm:2e-dual} and \ref{thm:scalar_prod}.

In general, spherical designs are multisets of points in the sphere. Designs of small cardinality such as tight or almost tight designs do not have coinciding points. However, it is interesting to determine the minimal size of a design that has two or, more generally, $\ell$ coinciding points. If a tight design exists, $\ell$ copies of this design may have $\ell$ coinciding points so a universal bound should not be above $\ell M(d,e)$. The following theorem shows that this lower bound follows from Theorem~\ref{thm:2e-dual}.

\begin{theorem}\label{thm:coinc}
Let $X$ be a spherical $2e$-design that contains $\ell$ copies of the same point. Then $|X|\geq\ell M(d,e)$.
\end{theorem}

\begin{proof}
Let $N=|X|$. For $\ell=1$, this is just the standard lower bound from \cite{DGS77}. For $\ell\geq 2$, we note that $X$ has at least $M(d,e)+1$ points so we can use Theorem~\ref{thm:2e-dual}. The existence of $\ell$ points whose pairwise inner products are 1 implies the existence of $\ell$ points in the dual $2$-design $Y$ whose inner products are $-\frac {M} {N-M}$ because $R_e(1)=M$. $Y$ contains a regular simplex with $\ell$ vertices with pairwise inner products $-\frac {M} {N-M}$ that cannot be smaller than $-\frac 1 {\ell-1}$. Therefore, $N-M\geq M(\ell-1)$ and $|X|\geq \ell M(d,e)$.
\end{proof}

It would be interesting to determine when the bound of Theorem~\ref{thm:coinc} is exact for $\ell\geq 2$. We note that the situation here is more complicated than for the case $\ell=1$ because there are no rigid constraints on inner products. Whenever a tight design exists any $\ell$ rotational copies of it attain equality in the bound, as long as they all share the same point.

In a similar manner, we can estimate the size of a $2e$-design if it has two antipodal points.

\begin{theorem}\label{thm:antipodal}
If $X$ is a spherical $2e$-design with two antipodal points, then $|X|\geq 2\binom{d+e-1}{d-1}$.
\end{theorem}

\begin{proof}
Let $x_1$ and $x_2$ be the two antipodal points in the design $X$ of size $N$ and $y_1$ and $y_2$ be their counterparts from $S^{N-M(d,e)-1}$ as described by Theorem~\ref{thm:2e-dual}. Note that a tight $2e$-design does not have antipodal points because $R_e(-1)\neq 0$ (see the calculation of $R_e(-1)$ below). Then
\[
\langle y_1, y_2 \rangle = - \frac {R_e(-1)} {N-M(d,e)}.
\]

As already mentioned, $R_e(x)=C_e^{(d/2)}(x)+C_{e-1}^{(d/2)}(x)$. Note that each polynomial $C_k^{(d/2)}$ is even for even $k$ and odd for odd $k$, so, using the values $C_k^{(d/2)}(1)$ from \cite[Theorem 2.4]{DGS77}, we get
\[
R_e(-1)=(-1)^e C_e^{(d/2)}(1)+(-1)^{e-1}C_{e-1}^{(d/2)}(1)=(-1)^e \binom{d+e-1}{d-1} + (-1)^{e-1} \binom{d+e-2}{d-1}.
\]
Therefore,
\[
N-M(d,e)\geq |R_e(-1)|=\binom{d+e-1}{d-1} - \binom{d+e-2}{d-1}.
\]
From the formula for $M(d,e)$ we get the required lower bound.
\end{proof}

Interestingly, the bound of Theorem~\ref{thm:antipodal} is the same as the Fisher-type bound for $(2e+1)$-designs from \cite{DGS77}. Tight $(2e+1)$-designs are necessarily antipodal so the bound is exact for them. It would also be interesting to determine when the bound of Theorem~\ref{thm:antipodal} is tight. We note that this bound is already present in the literature, see \cite[Corollary 4.5.11]{B02}.

In the case when the size of a design does not allow for two coinciding points as in Theorem~\ref{thm:coinc} or two antipodal points as in Theorem~\ref{thm:antipodal}, Theorem~\ref{thm:scalar_prod} provides an upper and lower bound, respectively, for possible inner products in a $2e$-design.

\begin{corollary}\label{cor:scalar_prod}
Let $M(d,e)\leq N \leq 2M(d,e)-1$. By $t_{d,e,N}$ we define the largest root of the equation $R_e(t)=N-M(d,e)$ for $t\in[-1,1]$. If $N<2\binom{d+e-1}{d-1}$, by $r_{d,e,N}$ we define the smallest root of the equation $R_e(r)=(-1)^e (N-M(d,e))$ for $r\in[-1,1]$. Otherwise, $r_{d,e,N}=-1$. Then all inner products of a spherical $2e$-design of size $N$ in dimension $d$ belong to the interval $[r_{d,e,N}, t_{d,e,N}]$.
\end{corollary}

\begin{proof}
This is an immediate consequence of Theorem~\ref{thm:scalar_prod} and the same calculations as in Theorem~\ref{thm:coinc} and Theorem~\ref{thm:antipodal}.
\end{proof}

Preliminary numerical investigations suggest that the bounds in Corollary~\ref{cor:scalar_prod} improve on some of the bounds in the literature \cite{BBD99, B02, BBKS09, BDHSS15}. We leave a detailed comparison to future work. In this paper, we show several cases in which these bounds are attained.

The upper bound in Corollary~\ref{cor:scalar_prod} is exact for tight $2e$-designs. The lower bound is exact for tight $(2e+1)$-designs as mentioned above. Both bounds are also exact for almost tight 2-designs by Mimura in even dimensions, whose inner products are precisely $0$ and $-\frac 2 d$. Below, we will show one more case, when the upper bound of Corollary~\ref{cor:scalar_prod} is exact.

\begin{theorem}\label{thm:5-design_mod}
For any 4-design of size $(d+1)d$ in dimension $d$, its inner products are no larger than $\frac d {d+2}$. If there exists a spherical 5-design of size $(d+1)d$ in dimension $d$, then this bound is exact.
\end{theorem}

For the proof of this theorem, we will need a lemma about orthogonal transformations within the latitude slices of a design. We note that it follows from the results in \cite{B92} and essentially uses the argument of \cite[Theorem 8.2]{DGS77} but the proof is very short and self-contained so it makes sense to include it here.

Let $X$ be a $\tau$-design in $S^{d-1}$. For a fixed $u\in S^{d-1}$, we split $X$ into the sets from latitudinal spheres, that is, for each $c\in [-1,1]$, we define $X_c=\{v\in X: \langle u,v \rangle = c\}$. For $c\in (-1,1)$, if $X_c$ is nonempty, we define $Z_c\subset S^{d-2}$ as the corresponding derived code such that $X_c=\{(c,\sqrt{1-c^2}v):v\in Z_c\}$. The lemma claims that if all derived codes are $\tau$-designs too, then we can rotate $X_c$ within their own latitudinal subspheres and preserve the $\tau$-design property of the whole set.

\begin{lemma}\label{lem:der_design}
Let $X$ be a $\tau$-design in $S^{d-1}$ and suppose that, for a fixed unit vector $u$, all nonempty derived codes $Z_c$, $-1<c<1$, are $\tau$-designs in $S^{d-2}$. Then for any orthogonal transformations $g_c\in O(d-1)$, the multiset
\[
X'=X_1\uplus X_{-1}\uplus \biguplus_{\substack{-1<c<1\\X_c\ne\varnothing}}X_c',
\]
where $X'_c=\{(c,\sqrt{1-c^2}v):v\in g_c(Z_c)\}$, is also a $\tau$-design.
\end{lemma}

\begin{proof}
Let $f$ be an arbitrary homogeneous polynomial of $d$ real variables of degree $k\leq \tau$. Then $f(t,\eta) = \sum_r t^{k-r} f_r(\eta)$, where $t$ stands for the coordinate with respect to $u$ and $\eta\in u^{\perp}$ and $f_r$ has degree $r$. Then
\[
\sum\limits_{v\in X'_c} f(v) = \sum_r c^{k-r} (1-c^2)^{r/2} \sum\limits_{v\in g_c (Z_c)} f_r(v) =  \sum_r c^{k-r} (1-c^2)^{r/2} \sum\limits_{v\in Z_c} f_r(v) = \sum\limits_{v\in X_c} f(v)
\]
Here we use the fact that both $Z_c$ and $g_c(Z_c)$ are $\tau$-designs in $S^{d-2}$ so $\sum f_r(v)$ is the same for both of them. Summing these equalities over $c$, we get that $X'$ is also a $\tau$-design. The north and south poles, if present in $X$, are unchanged and contribute equally to both sums.
\end{proof}

\begin{proof}[Proof of Theorem~\ref{thm:5-design_mod}]
The upper bound $\frac d {d+2}$ is obtained from Corollary~\ref{cor:scalar_prod}. Now we assume there exists a tight 5-design in dimension $d$. A tight 5-design is antipodal with inner products $-1, \pm \frac 1 {\sqrt{d+2}}$ \cite[Theorem 5.12]{DGS77}. Fix a point $u$ of the tight 5-design as the north pole. Derived codes in latitudes $\frac 1 {\sqrt{d+2}}$ and $-\frac 1 {\sqrt{d+2}}$ are necessarily 4-designs \cite[Theorem 8.2]{DGS77} (moreover, just counting their size we can see they are tight 4-designs). By Lemma~\ref{lem:der_design}, rotating the derived code in one of the latitudinal subspheres we still get a 4-design. We consider an arbitrary orthogonal transformation making two of the points align, that is, the north pole, the south pole, a point with latitude $\frac 1 {\sqrt{d+2}}$, and a point with latitude $-\frac 1 {\sqrt{d+2}}$ are in the same 2-dimensional subspace and the points with latitudes $\pm\frac 1 {\sqrt{d+2}}$ are not antipodal. Then the angle between the points with latitudes $\pm\frac 1 {\sqrt{d+2}}$ is $\pi-2\alpha$, where $\cos\alpha=\frac 1 {\sqrt{d+2}}$ and, therefore, their inner product is $1-2(\cos\alpha)^2=\frac d {d+2}$, as required.
\end{proof}

Tight 5-designs in dimension $d$ are known to exist if and only if tight 4-designs exist in dimension $d-1$ \cite{DGS77}. In particular, they exist in dimensions 3, 7, and 23, which provide exact bounds for Corollary~\ref{cor:scalar_prod}. The case of dimension 3 is particularly interesting. As we showed above, the smallest 4-designs have size 12 in $S^2$. Hardin and Sloane noticed that configurations obtained by rotations of one of the pentagons are also 4-designs \cite{HS92}. Due to Theorem~\ref{thm:5-design_mod}, the only configuration among these with the maximal inner product $3/5$ is the one with the north pole, the south pole, and two aligned pentagons. In dimensions 7 and 23, there are infinitely many optimal 4-designs of that type because aligning two points in codimension 1 does not guarantee that all other points are aligned as well.

\bigskip

\noindent 
\textbf{Note on Related Work.}
After obtaining the main results of this paper, we became aware of
the preprint~\cite{SS26}, which also studies spherical designs with
cardinalities close to the Fisher-type lower bound and discusses
related dual constructions. We mention this work for completeness.
Our proofs were developed independently and do not rely on its claims.

\bigskip

\noindent
\textbf{Acknowledgments.} 
This research was initiated during A. Glazyrin's visit to Aichi University of Education, which was supported by the JSPS--NSF Supplemental Research Collaboration Opportunity in Japan (JPJSSRC20260008). H. Nozaki was partially supported by JSPS KAKENHI Grant Numbers JP22K03402, JP24K06688, and JP26K00617. 
A. Glazyrin was partially supported by the National Science Foundation grant DMS-2349063, the NSF-JSPS supplementary grant covering travel expenses to Japan, and the Simons Foundation’s Travel Support for Mathematicians program. 

 \bigskip

 \noindent 
 \textbf{Data Availability Statement.}
All computational data supporting the results of this paper, together with the programs that verify them, are publicly available at \url{https://github.com/alexey-glazyrin/design_verifier}. 

\bigskip

\noindent
\textbf{Conflict of Interest.} 
The authors state that there is no conflict of interest. 

\bigskip

\noindent
\textbf{Declaration of Generative AI Use.}
Generative AI tools (OpenAI ChatGPT Sol and Astra, Claude Fable) were used for exploratory research discussions, searching the literature, and editing the text. AI assistance was indispensable for writing and executing semidefinite programs that led to the computational results. Claude Fable, through Claude Code, was used by the authors to assist with creating the code for the semidefinite programs in Python, for the rounding and verification procedures, and for running the code.
All main mathematical statements and proofs are due to the authors. The statements and calculations produced by AI were verified by the authors, who take full responsibility for their correctness.

\begin{appendices}
\section{Table of certificates of nonexistence}\label{app:certificates}

\begin{table}[!ht]\centering\footnotesize
\caption{Exact SDP certificates: for every row there is no spherical $2e$-design with $N=M(d,e)+s$ points in $S^{d-1}$. $D$ is the kernel truncation parameter. $B$ denotes the exact bound of the certificate; the displayed values of $B$ and $N-B$ are rounded to nine decimal places.} \label{tab:certificates}
\vspace{6pt}
\begin{tabular}{@{}rrrrrrrr@{}}
\toprule
$e$ & $d$ & $s$ & $M(d,e)$ & $N$ & $D$ & \multicolumn{1}{c}{$B$} & \multicolumn{1}{c@{}}{$N-B$} \\
\midrule
2 & 3 & 2 & 9 & 11 & 6 & 10.990000000 & 0.010000000 \\
\addlinespace[2pt]
2 & 4 & 2 & 14 & 16 & 12 & 15.800000000 & 0.200000000 \\
\addlinespace[2pt]
2 & 5 & 2 & 20 & 22 & 8 & 21.800000000 & 0.200000000 \\
\midrule
3 & 3 & 2 & 16 & 18 & 8 & 17.970000000 & 0.030000000 \\
\addlinespace[2pt]
3 & 4 & 3 & 30 & 33 & 10 & 32.800000000 & 0.200000000 \\
\addlinespace[2pt]
3 & 5 & 4 & 50 & 54 & 10 & 53.800000000 & 0.200000000 \\
\addlinespace[2pt]
3 & 10 & 8 & 275 & 283 & 8 & 282.949999999 & 0.050000001 \\
3 & 10 & 9 & 275 & 284 & 8 & 283.949999990 & 0.050000010 \\
\addlinespace[2pt]
3 & 11 & 2 & 352 & 354 & 6 & 353.800000000 & 0.200000000 \\
3 & 11 & 3 & 352 & 355 & 6 & 353.999999872 & 1.000000128 \\
3 & 11 & 4 & 352 & 356 & 6 & 354.999999993 & 1.000000007 \\
3 & 11 & 5 & 352 & 357 & 6 & 356.000000000 & 1.000000000 \\
3 & 11 & 6 & 352 & 358 & 8 & 356.999998006 & 1.000001994 \\
3 & 11 & 7 & 352 & 359 & 8 & 357.999985971 & 1.000014029 \\
3 & 11 & 8 & 352 & 360 & 8 & 358.999995700 & 1.000004300 \\
3 & 11 & 9 & 352 & 361 & 8 & 359.999956206 & 1.000043794 \\
3 & 11 & 10 & 352 & 362 & 8 & 360.999990603 & 1.000009397 \\
\addlinespace[2pt]
3 & 12 & 2 & 442 & 444 & 6 & 443.800000000 & 0.200000000 \\
3 & 12 & 3 & 442 & 445 & 6 & 443.999999747 & 1.000000253 \\
3 & 12 & 4 & 442 & 446 & 6 & 444.999999995 & 1.000000005 \\
3 & 12 & 5 & 442 & 447 & 6 & 445.999999998 & 1.000000002 \\
3 & 12 & 6 & 442 & 448 & 6 & 446.999999995 & 1.000000005 \\
3 & 12 & 7 & 442 & 449 & 8 & 447.999994216 & 1.000005784 \\
3 & 12 & 8 & 442 & 450 & 8 & 448.999975014 & 1.000024986 \\
3 & 12 & 9 & 442 & 451 & 8 & 449.999960641 & 1.000039359 \\
3 & 12 & 10 & 442 & 452 & 8 & 450.999795618 & 1.000204382 \\
3 & 12 & 11 & 442 & 453 & 8 & 452.799265865 & 0.200734135 \\
\bottomrule
\end{tabular}
\end{table}

\begin{table}[p]\centering\footnotesize
\addtocounter{table}{-1}
\caption{Continued (part 2 of 6).}
\vspace{6pt}
\begin{tabular}{@{}rrrrrrrr@{}}
\toprule
$e$ & $d$ & $s$ & $M(d,e)$ & $N$ & $D$ & \multicolumn{1}{c}{$B$} & \multicolumn{1}{c@{}}{$N-B$} \\
\midrule
3 & 13 & 2 & 546 & 548 & 6 & 547.800000000 & 0.200000000 \\
3 & 13 & 3 & 546 & 549 & 6 & 547.999999755 & 1.000000245 \\
3 & 13 & 4 & 546 & 550 & 6 & 548.999999973 & 1.000000027 \\
3 & 13 & 5 & 546 & 551 & 6 & 549.999999990 & 1.000000010 \\
3 & 13 & 6 & 546 & 552 & 6 & 551.499999968 & 0.500000032 \\
3 & 13 & 7 & 546 & 553 & 8 & 551.999994520 & 1.000005480 \\
3 & 13 & 8 & 546 & 554 & 8 & 552.999978669 & 1.000021331 \\
3 & 13 & 9 & 546 & 555 & 8 & 553.999976404 & 1.000023596 \\
3 & 13 & 10 & 546 & 556 & 8 & 554.999995284 & 1.000004716 \\
3 & 13 & 11 & 546 & 557 & 8 & 555.999937877 & 1.000062123 \\
\addlinespace[2pt]
3 & 14 & 2 & 665 & 667 & 6 & 666.800000000 & 0.200000000 \\
3 & 14 & 3 & 665 & 668 & 6 & 666.988543587 & 1.011456413 \\
3 & 14 & 4 & 665 & 669 & 6 & 667.999999967 & 1.000000033 \\
3 & 14 & 5 & 665 & 670 & 6 & 668.999999955 & 1.000000045 \\
3 & 14 & 6 & 665 & 671 & 8 & 669.999965005 & 1.000034995 \\
3 & 14 & 7 & 665 & 672 & 8 & 670.999987441 & 1.000012559 \\
3 & 14 & 8 & 665 & 673 & 8 & 671.999881379 & 1.000118621 \\
3 & 14 & 9 & 665 & 674 & 8 & 672.999936889 & 1.000063111 \\
3 & 14 & 10 & 665 & 675 & 8 & 673.999990743 & 1.000009257 \\
3 & 14 & 11 & 665 & 676 & 8 & 675.799942123 & 0.200057877 \\
\addlinespace[2pt]
3 & 15 & 2 & 800 & 802 & 6 & 801.489684553 & 0.510315447 \\
3 & 15 & 3 & 800 & 803 & 6 & 802.499892390 & 0.500107610 \\
3 & 15 & 4 & 800 & 804 & 6 & 803.499957056 & 0.500042944 \\
3 & 15 & 5 & 800 & 805 & 6 & 804.499973977 & 0.500026023 \\
3 & 15 & 6 & 800 & 806 & 8 & 804.999932820 & 1.000067180 \\
3 & 15 & 7 & 800 & 807 & 8 & 805.999921027 & 1.000078973 \\
3 & 15 & 8 & 800 & 808 & 8 & 806.999976087 & 1.000023913 \\
3 & 15 & 9 & 800 & 809 & 8 & 807.999977694 & 1.000022306 \\
3 & 15 & 10 & 800 & 810 & 8 & 809.799873485 & 0.200126515 \\
\addlinespace[2pt]
3 & 16 & 2 & 952 & 954 & 6 & 953.498055512 & 0.501944488 \\
3 & 16 & 3 & 952 & 955 & 6 & 953.999965002 & 1.000034998 \\
3 & 16 & 4 & 952 & 956 & 8 & 954.999907180 & 1.000092820 \\
3 & 16 & 5 & 952 & 957 & 8 & 955.999446554 & 1.000553446 \\
3 & 16 & 6 & 952 & 958 & 8 & 956.999878604 & 1.000121396 \\
3 & 16 & 7 & 952 & 959 & 8 & 958.798997901 & 0.201002099 \\
\midrule
4 & 3 & 2 & 25 & 27 & 10 & 26.979999999 & 0.020000001 \\
\addlinespace[2pt]
4 & 18 & 2 & 7125 & 7127 & 8 & 7061.813998027 & 65.186001973 \\
4 & 18 & 3 & 7125 & 7128 & 8 & 7127.752437914 & 0.247562086 \\
4 & 18 & 4 & 7125 & 7129 & 8 & 7127.974258143 & 1.025741857 \\
4 & 18 & 5 & 7125 & 7130 & 8 & 7128.908821029 & 1.091178971 \\
4 & 18 & 6 & 7125 & 7131 & 8 & 7129.445068064 & 1.554931936 \\
4 & 18 & 7 & 7125 & 7132 & 8 & 7131.735870455 & 0.264129545 \\
4 & 18 & 8 & 7125 & 7133 & 8 & 7131.970523951 & 1.029476049 \\
4 & 18 & 9 & 7125 & 7134 & 8 & 7132.938355976 & 1.061644024 \\
4 & 18 & 10 & 7125 & 7135 & 8 & 7133.911366948 & 1.088633052 \\
4 & 18 & 11 & 7125 & 7136 & 8 & 7134.978957268 & 1.021042732 \\
4 & 18 & 12 & 7125 & 7137 & 8 & 7135.998062577 & 1.001937423 \\
4 & 18 & 13 & 7125 & 7138 & 8 & 7136.999971146 & 1.000028854 \\
\bottomrule
\end{tabular}
\end{table}

\begin{table}[p]\centering\footnotesize
\addtocounter{table}{-1}
\caption{Continued (part 3 of 6).}
\vspace{6pt}
\begin{tabular}{@{}rrrrrrrr@{}}
\toprule
$e$ & $d$ & $s$ & $M(d,e)$ & $N$ & $D$ & \multicolumn{1}{c}{$B$} & \multicolumn{1}{c@{}}{$N-B$} \\
\midrule
4 & 18 & 14 & 7125 & 7139 & 8 & 7137.996795719 & 1.003204281 \\
4 & 18 & 15 & 7125 & 7140 & 8 & 7138.999447263 & 1.000552737 \\
4 & 18 & 16 & 7125 & 7141 & 8 & 7139.995864754 & 1.004135246 \\
4 & 18 & 17 & 7125 & 7142 & 8 & 7140.916782370 & 1.083217630 \\
4 & 18 & 18 & 7125 & 7143 & 8 & 7141.986682674 & 1.013317326 \\
4 & 18 & 19 & 7125 & 7144 & 8 & 7142.946258378 & 1.053741622 \\
4 & 18 & 20 & 7125 & 7145 & 8 & 7143.995507216 & 1.004492784 \\
4 & 18 & 21 & 7125 & 7146 & 8 & 7144.974426898 & 1.025573102 \\
4 & 18 & 22 & 7125 & 7147 & 8 & 7145.969741048 & 1.030258952 \\
4 & 18 & 23 & 7125 & 7148 & 8 & 7146.858709040 & 1.141290960 \\
\addlinespace[2pt]
4 & 19 & 2 & 8645 & 8647 & 8 & 8645.431390380 & 1.568609620 \\
4 & 19 & 3 & 8645 & 8648 & 8 & 8646.996632020 & 1.003367980 \\
4 & 19 & 4 & 8645 & 8649 & 8 & 8647.987463708 & 1.012536292 \\
4 & 19 & 5 & 8645 & 8650 & 8 & 8648.122196088 & 1.877803912 \\
4 & 19 & 6 & 8645 & 8651 & 8 & 8649.986182235 & 1.013817765 \\
4 & 19 & 7 & 8645 & 8652 & 8 & 8651.757578456 & 0.242421544 \\
4 & 19 & 8 & 8645 & 8653 & 8 & 8651.964067789 & 1.035932211 \\
4 & 19 & 9 & 8645 & 8654 & 8 & 8652.823791702 & 1.176208298 \\
4 & 19 & 10 & 8645 & 8655 & 8 & 8653.919230133 & 1.080769867 \\
4 & 19 & 11 & 8645 & 8656 & 8 & 8654.927632299 & 1.072367701 \\
4 & 19 & 12 & 8645 & 8657 & 8 & 8655.843173770 & 1.156826230 \\
4 & 19 & 13 & 8645 & 8658 & 8 & 8656.926957959 & 1.073042041 \\
4 & 19 & 14 & 8645 & 8659 & 8 & 8657.988991850 & 1.011008150 \\
4 & 19 & 15 & 8645 & 8660 & 8 & 8658.919771546 & 1.080228454 \\
4 & 19 & 16 & 8645 & 8661 & 8 & 8659.845876515 & 1.154123485 \\
4 & 19 & 17 & 8645 & 8662 & 8 & 8660.999099598 & 1.000900402 \\
4 & 19 & 18 & 8645 & 8663 & 8 & 8661.996380593 & 1.003619407 \\
4 & 19 & 19 & 8645 & 8664 & 8 & 8662.974630599 & 1.025369401 \\
4 & 19 & 20 & 8645 & 8665 & 8 & 8663.946896470 & 1.053103530 \\
4 & 19 & 21 & 8645 & 8666 & 8 & 8664.972493515 & 1.027506485 \\
4 & 19 & 22 & 8645 & 8667 & 8 & 8665.998933804 & 1.001066196 \\
4 & 19 & 23 & 8645 & 8668 & 8 & 8666.997531804 & 1.002468196 \\
\addlinespace[2pt]
4 & 20 & 2 & 10395 & 10397 & 8 & 10395.508447901 & 1.491552099 \\
4 & 20 & 3 & 10395 & 10398 & 8 & 10396.842571659 & 1.157428341 \\
4 & 20 & 4 & 10395 & 10399 & 8 & 10397.943984256 & 1.056015744 \\
4 & 20 & 5 & 10395 & 10400 & 8 & 10398.991998066 & 1.008001934 \\
4 & 20 & 6 & 10395 & 10401 & 8 & 10399.995544281 & 1.004455719 \\
4 & 20 & 7 & 10395 & 10402 & 8 & 10400.996772428 & 1.003227572 \\
4 & 20 & 8 & 10395 & 10403 & 8 & 10401.602189592 & 1.397810408 \\
4 & 20 & 9 & 10395 & 10404 & 8 & 10402.994700564 & 1.005299436 \\
4 & 20 & 10 & 10395 & 10405 & 8 & 10402.322348271 & 2.677651729 \\
4 & 20 & 11 & 10395 & 10406 & 8 & 10404.998110825 & 1.001889175 \\
4 & 20 & 12 & 10395 & 10407 & 8 & 10405.695714536 & 1.304285464 \\
4 & 20 & 13 & 10395 & 10408 & 8 & 10406.999017200 & 1.000982800 \\
4 & 20 & 14 & 10395 & 10409 & 8 & 10407.798130454 & 1.201869546 \\
4 & 20 & 15 & 10395 & 10410 & 8 & 10408.894224121 & 1.105775879 \\
4 & 20 & 16 & 10395 & 10411 & 8 & 10409.485217613 & 1.514782387 \\
4 & 20 & 17 & 10395 & 10412 & 8 & 10410.961710211 & 1.038289789 \\
4 & 20 & 18 & 10395 & 10413 & 8 & 10411.993228450 & 1.006771550 \\
\bottomrule
\end{tabular}
\end{table}

\begin{table}[p]\centering\footnotesize
\addtocounter{table}{-1}
\caption{Continued (part 4 of 6).}
\vspace{6pt}
\begin{tabular}{@{}rrrrrrrr@{}}
\toprule
$e$ & $d$ & $s$ & $M(d,e)$ & $N$ & $D$ & \multicolumn{1}{c}{$B$} & \multicolumn{1}{c@{}}{$N-B$} \\
\midrule
4 & 20 & 19 & 10395 & 10414 & 8 & 10412.854139347 & 1.145860653 \\
4 & 20 & 20 & 10395 & 10415 & 8 & 10413.961964949 & 1.038035051 \\
4 & 20 & 21 & 10395 & 10416 & 8 & 10414.999653824 & 1.000346176 \\
4 & 20 & 22 & 10395 & 10417 & 8 & 10415.958467273 & 1.041532727 \\
\addlinespace[2pt]
4 & 21 & 2 & 12397 & 12399 & 8 & 12396.736440739 & 2.263559261 \\
4 & 21 & 3 & 12397 & 12400 & 8 & 12398.921977337 & 1.078022663 \\
4 & 21 & 4 & 12397 & 12401 & 8 & 12399.022321976 & 1.977678024 \\
4 & 21 & 5 & 12397 & 12402 & 8 & 12400.961241791 & 1.038758209 \\
4 & 21 & 6 & 12397 & 12403 & 8 & 12401.821754871 & 1.178245129 \\
4 & 21 & 7 & 12397 & 12404 & 8 & 12402.992679529 & 1.007320471 \\
4 & 21 & 8 & 12397 & 12405 & 8 & 12403.996455941 & 1.003544059 \\
4 & 21 & 9 & 12397 & 12406 & 8 & 12404.866595498 & 1.133404502 \\
4 & 21 & 10 & 12397 & 12407 & 8 & 12405.418256653 & 1.581743347 \\
4 & 21 & 11 & 12397 & 12408 & 8 & 12406.647382369 & 1.352617631 \\
4 & 21 & 12 & 12397 & 12409 & 8 & 12407.643758752 & 1.356241248 \\
4 & 21 & 13 & 12397 & 12410 & 8 & 12408.845340709 & 1.154659291 \\
4 & 21 & 14 & 12397 & 12411 & 8 & 12409.498584953 & 1.501415047 \\
4 & 21 & 15 & 12397 & 12412 & 8 & 12410.992477650 & 1.007522350 \\
4 & 21 & 16 & 12397 & 12413 & 8 & 12411.997676468 & 1.002323532 \\
4 & 21 & 17 & 12397 & 12414 & 8 & 12412.967167616 & 1.032832384 \\
4 & 21 & 18 & 12397 & 12415 & 8 & 12413.938630050 & 1.061369950 \\
4 & 21 & 19 & 12397 & 12416 & 8 & 12414.872663574 & 1.127336426 \\
4 & 21 & 20 & 12397 & 12417 & 8 & 12415.694215313 & 1.305784687 \\
4 & 21 & 21 & 12397 & 12418 & 8 & 12416.999346772 & 1.000653228 \\
4 & 21 & 22 & 12397 & 12419 & 8 & 12417.996724238 & 1.003275762 \\
4 & 21 & 23 & 12397 & 12420 & 8 & 12418.995283679 & 1.004716321 \\
\addlinespace[2pt]
4 & 22 & 2 & 14674 & 14676 & 8 & 14675.364631164 & 0.635368836 \\
4 & 22 & 3 & 14674 & 14677 & 8 & 14675.944035758 & 1.055964242 \\
4 & 22 & 4 & 14674 & 14678 & 8 & 14676.979018213 & 1.020981787 \\
4 & 22 & 5 & 14674 & 14679 & 8 & 14677.836409596 & 1.163590404 \\
4 & 22 & 6 & 14674 & 14680 & 8 & 14678.966270174 & 1.033729826 \\
4 & 22 & 7 & 14674 & 14681 & 8 & 14679.065422132 & 1.934577868 \\
4 & 22 & 8 & 14674 & 14682 & 8 & 14680.995934839 & 1.004065161 \\
4 & 22 & 9 & 14674 & 14683 & 8 & 14681.892782367 & 1.107217633 \\
4 & 22 & 10 & 14674 & 14684 & 8 & 14682.770763541 & 1.229236459 \\
4 & 22 & 11 & 14674 & 14685 & 8 & 14683.919953615 & 1.080046385 \\
4 & 22 & 12 & 14674 & 14686 & 8 & 14684.245182456 & 1.754817544 \\
4 & 22 & 13 & 14674 & 14687 & 8 & 14685.347458420 & 1.652541580 \\
4 & 22 & 14 & 14674 & 14688 & 8 & 14686.764992103 & 1.235007897 \\
4 & 22 & 15 & 14674 & 14689 & 8 & 14687.020549913 & 1.979450087 \\
4 & 22 & 16 & 14674 & 14690 & 8 & 14687.562477250 & 2.437522750 \\
4 & 22 & 17 & 14674 & 14691 & 8 & 14689.530434639 & 1.469565361 \\
4 & 22 & 18 & 14674 & 14692 & 8 & 14690.573261602 & 1.426738398 \\
4 & 22 & 19 & 14674 & 14693 & 8 & 14690.683789069 & 2.316210931 \\
4 & 22 & 20 & 14674 & 14694 & 8 & 14692.876278267 & 1.123721733 \\
4 & 22 & 21 & 14674 & 14695 & 8 & 14693.998446855 & 1.001553145 \\
\bottomrule
\end{tabular}
\end{table}

\begin{table}[p]\centering\footnotesize
\addtocounter{table}{-1}
\caption{Continued (part 5 of 6).}
\vspace{6pt}
\begin{tabular}{@{}rrrrrrrr@{}}
\toprule
$e$ & $d$ & $s$ & $M(d,e)$ & $N$ & $D$ & \multicolumn{1}{c}{$B$} & \multicolumn{1}{c@{}}{$N-B$} \\
\midrule
4 & 23 & 2 & 17250 & 17252 & 8 & 17251.800000010 & 0.199999990 \\
4 & 23 & 3 & 17250 & 17253 & 8 & 17251.628926270 & 1.371073730 \\
4 & 23 & 4 & 17250 & 17254 & 8 & 17252.733492188 & 1.266507812 \\
4 & 23 & 5 & 17250 & 17255 & 8 & 17253.794609713 & 1.205390287 \\
4 & 23 & 6 & 17250 & 17256 & 8 & 17254.529285017 & 1.470714983 \\
4 & 23 & 7 & 17250 & 17257 & 8 & 17255.825157963 & 1.174842037 \\
4 & 23 & 8 & 17250 & 17258 & 8 & 17256.969975915 & 1.030024085 \\
4 & 23 & 9 & 17250 & 17259 & 8 & 17257.290644997 & 1.709355003 \\
4 & 23 & 10 & 17250 & 17260 & 8 & 17258.895930486 & 1.104069514 \\
4 & 23 & 11 & 17250 & 17261 & 8 & 17259.977936481 & 1.022063519 \\
4 & 23 & 12 & 17250 & 17262 & 8 & 17260.870831500 & 1.129168500 \\
4 & 23 & 13 & 17250 & 17263 & 8 & 17261.983742720 & 1.016257280 \\
4 & 23 & 14 & 17250 & 17264 & 8 & 17262.676669139 & 1.323330861 \\
4 & 23 & 15 & 17250 & 17265 & 8 & 17263.996280108 & 1.003719892 \\
4 & 23 & 16 & 17250 & 17266 & 8 & 17264.855403790 & 1.144596210 \\
4 & 23 & 17 & 17250 & 17267 & 8 & 17265.096704276 & 1.903295724 \\
4 & 23 & 18 & 17250 & 17268 & 8 & 17266.985778009 & 1.014221991 \\
4 & 23 & 19 & 17250 & 17269 & 8 & 17266.627026374 & 2.372973626 \\
4 & 23 & 20 & 17250 & 17270 & 8 & 17268.987541604 & 1.012458396 \\
\addlinespace[2pt]
4 & 24 & 2 & 20150 & 20152 & 8 & 20151.799999988 & 0.200000012 \\
4 & 24 & 3 & 20150 & 20153 & 8 & 20151.561690544 & 1.438309456 \\
4 & 24 & 4 & 20150 & 20154 & 8 & 20152.738458242 & 1.261541758 \\
4 & 24 & 5 & 20150 & 20155 & 8 & 20153.904902873 & 1.095097127 \\
4 & 24 & 6 & 20150 & 20156 & 8 & 20150.389380220 & 5.610619780 \\
4 & 24 & 7 & 20150 & 20157 & 8 & 20155.944531460 & 1.055468540 \\
4 & 24 & 8 & 20150 & 20158 & 8 & 20156.750274651 & 1.249725349 \\
4 & 24 & 9 & 20150 & 20159 & 8 & 20156.617691976 & 2.382308024 \\
4 & 24 & 10 & 20150 & 20160 & 8 & 20156.248599387 & 3.751400613 \\
4 & 24 & 11 & 20150 & 20161 & 8 & 20156.962927698 & 4.037072302 \\
4 & 24 & 12 & 20150 & 20162 & 8 & 20160.884959614 & 1.115040386 \\
4 & 24 & 13 & 20150 & 20163 & 8 & 20157.617259582 & 5.382740418 \\
4 & 24 & 14 & 20150 & 20164 & 8 & 20162.867972203 & 1.132027797 \\
4 & 24 & 15 & 20150 & 20165 & 8 & 20163.969471007 & 1.030528993 \\
4 & 24 & 16 & 20150 & 20166 & 8 & 20161.080148389 & 4.919851611 \\
4 & 24 & 17 & 20150 & 20167 & 8 & 20165.879651576 & 1.120348424 \\
4 & 24 & 18 & 20150 & 20168 & 8 & 20166.614501429 & 1.385498571 \\
4 & 24 & 19 & 20150 & 20169 & 8 & 20168.783487323 & 0.216512677 \\
\addlinespace[2pt]
4 & 25 & 2 & 23400 & 23402 & 8 & 23401.800000003 & 0.199999997 \\
4 & 25 & 3 & 23400 & 23403 & 8 & 23402.134840712 & 0.865159288 \\
4 & 25 & 4 & 23400 & 23404 & 8 & 23401.822876892 & 2.177123108 \\
4 & 25 & 5 & 23400 & 23405 & 8 & 23403.483709245 & 1.516290755 \\
4 & 25 & 6 & 23400 & 23406 & 8 & 23403.628362344 & 2.371637656 \\
4 & 25 & 7 & 23400 & 23407 & 8 & 23405.924229072 & 1.075770928 \\
4 & 25 & 8 & 23400 & 23408 & 8 & 23406.974034789 & 1.025965211 \\
4 & 25 & 9 & 23400 & 23409 & 8 & 23407.305226535 & 1.694773465 \\
4 & 25 & 10 & 23400 & 23410 & 8 & 23407.731487941 & 2.268512059 \\
4 & 25 & 11 & 23400 & 23411 & 8 & 23406.322327455 & 4.677672545 \\
4 & 25 & 12 & 23400 & 23412 & 8 & 23410.754877534 & 1.245122466 \\
4 & 25 & 13 & 23400 & 23413 & 8 & 23411.348503185 & 1.651496815 \\
\bottomrule
\end{tabular}
\end{table}

\begin{table}[p]\centering\footnotesize
\addtocounter{table}{-1}
\caption{Continued (part 6 of 6).}
\vspace{6pt}
\begin{tabular}{@{}rrrrrrrr@{}}
\toprule
$e$ & $d$ & $s$ & $M(d,e)$ & $N$ & $D$ & \multicolumn{1}{c}{$B$} & \multicolumn{1}{c@{}}{$N-B$} \\
\midrule
4 & 25 & 14 & 23400 & 23414 & 8 & 23412.985899215 & 1.014100785 \\
4 & 25 & 15 & 23400 & 23415 & 8 & 23413.957332610 & 1.042667390 \\
4 & 25 & 16 & 23400 & 23416 & 8 & 23412.837825605 & 3.162174395 \\
4 & 25 & 17 & 23400 & 23417 & 8 & 23413.244013641 & 3.755986359 \\
\addlinespace[2pt]
4 & 26 & 2 & 27027 & 27029 & 8 & 27028.799999925 & 0.200000075 \\
4 & 26 & 3 & 27027 & 27030 & 8 & 27029.799999817 & 0.200000183 \\
4 & 26 & 4 & 27027 & 27031 & 8 & 27030.800000251 & 0.199999749 \\
4 & 26 & 5 & 27027 & 27032 & 8 & 27031.800000784 & 0.199999216 \\
4 & 26 & 6 & 27027 & 27033 & 8 & 27031.273445163 & 1.726554837 \\
4 & 26 & 7 & 27027 & 27034 & 8 & 27031.298275014 & 2.701724986 \\
4 & 26 & 8 & 27027 & 27035 & 8 & 27033.876079131 & 1.123920869 \\
4 & 26 & 9 & 27027 & 27036 & 8 & 27019.965287208 & 16.034712792 \\
4 & 26 & 10 & 27027 & 27037 & 8 & 27035.313177535 & 1.686822465 \\
4 & 26 & 11 & 27027 & 27038 & 8 & 27027.449319223 & 10.550680777 \\
4 & 26 & 12 & 27027 & 27039 & 8 & 27037.883126620 & 1.116873380 \\
\midrule
5 & 27 & 2 & 197316 & 197318 & 10 & 197317.799999987 & 0.200000013 \\
\addlinespace[2pt]
5 & 28 & 2 & 232841 & 232843 & 10 & 232842.799999991 & 0.200000009 \\
5 & 28 & 3 & 232841 & 232844 & 10 & 232843.800000000 & 0.200000000 \\
5 & 28 & 4 & 232841 & 232845 & 10 & 232844.800000000 & 0.200000000 \\
5 & 28 & 5 & 232841 & 232846 & 10 & 232845.799998344 & 0.200001656 \\
5 & 28 & 6 & 232841 & 232847 & 10 & 232846.800086421 & 0.199913579 \\
5 & 28 & 7 & 232841 & 232848 & 10 & 232847.800043303 & 0.199956697 \\
5 & 28 & 8 & 232841 & 232849 & 10 & 232848.799614566 & 0.200385434 \\
5 & 28 & 9 & 232841 & 232850 & 10 & 232849.800036459 & 0.199963541 \\
5 & 28 & 10 & 232841 & 232851 & 10 & 232850.799897542 & 0.200102458 \\
5 & 28 & 11 & 232841 & 232852 & 10 & 232851.799975453 & 0.200024547 \\
5 & 28 & 12 & 232841 & 232853 & 10 & 232852.800104367 & 0.199895633 \\
5 & 28 & 13 & 232841 & 232854 & 10 & 232853.800019105 & 0.199980895 \\
5 & 28 & 14 & 232841 & 232855 & 10 & 232854.799995842 & 0.200004158 \\
5 & 28 & 15 & 232841 & 232856 & 10 & 232855.800010922 & 0.199989078 \\
5 & 28 & 16 & 232841 & 232857 & 10 & 232856.799983553 & 0.200016447 \\
5 & 28 & 17 & 232841 & 232858 & 10 & 232857.799985666 & 0.200014334 \\
5 & 28 & 18 & 232841 & 232859 & 10 & 232858.800013553 & 0.199986447 \\
5 & 28 & 19 & 232841 & 232860 & 10 & 232859.800031107 & 0.199968893 \\
5 & 28 & 20 & 232841 & 232861 & 10 & 232860.800065252 & 0.199934748 \\
5 & 28 & 21 & 232841 & 232862 & 10 & 232861.800415563 & 0.199584437 \\
5 & 28 & 22 & 232841 & 232863 & 10 & 232862.799626232 & 0.200373768 \\
5 & 28 & 23 & 232841 & 232864 & 10 & 232863.800050489 & 0.199949511 \\
5 & 28 & 24 & 232841 & 232865 & 10 & 232864.800716661 & 0.199283339 \\
5 & 28 & 25 & 232841 & 232866 & 10 & 232865.800306267 & 0.199693733 \\
5 & 28 & 26 & 232841 & 232867 & 10 & 232866.799896900 & 0.200103100 \\
5 & 28 & 27 & 232841 & 232868 & 10 & 232867.800451933 & 0.199548067 \\
5 & 28 & 28 & 232841 & 232869 & 10 & 232868.800557555 & 0.199442445 \\
5 & 28 & 29 & 232841 & 232870 & 10 & 232869.800479496 & 0.199520504 \\
5 & 28 & 30 & 232841 & 232871 & 10 & 232870.800145551 & 0.199854449 \\
\addlinespace[2pt]
5 & 29 & 2 & 273296 & 273298 & 10 & 273297.800000000 & 0.200000000 \\
\bottomrule
\end{tabular}
\end{table}
\newpage
\section{Linear programming bounds used in Table~\ref{tab:bounds}}

Three of the bounds of Table~\ref{tab:bounds} rely on the linear programming bound of~\cite{DGS77} for the smallest sizes. To keep Table~\ref{tab:bounds} self-contained, we record polynomials that prove what is needed. Recall from Subsection~\ref{subsect:certs} that if $f=\sum_ka_kP^d_k$ satisfies $a_0>0$, $a_k\le 0$ for $k>2e$, and $f(t)\ge 0$ on $[-1,1]$, then every spherical $2e$-design in $S^{d-1}$ has at least $f(1)/a_0$ points. Table~\ref{tab:lp} lists three such polynomials with integer coefficients, of degree $11$, $9$ and $9$. Each of them is positive on $[-1,1]$; this is verified exactly with a Sturm sequence. Since $f(1)/a_0$ exceeds $32$, $53$ and $282$, respectively, every $6$-design has at least $33$ points in $S^3$, at least $54$ points in $S^4$, and at least $283$ points in $S^9$. With polynomials of degree up to $40$ we found no linear programming bound that excludes the next size in any of the three cases.

\begin{table}[ht!]\centering
\caption{Polynomials $f=\sum_ka_kP^d_k$ for the linear programming bound $N\ge f(1)/a_0$ for spherical $6$-designs in $S^{d-1}$. The coefficients $a_7$ and $a_8$ are zero.}\label{tab:lp}
\vspace{6pt}
\begin{tabular}{@{}lrrr@{}}
\toprule
 & $d=4$ & $d=5$ & $d=10$ \\
\midrule
$a_0$ & $117935$ & $70$ & $104976$ \\
$a_1$ & $364840$ & $245$ & $621687$ \\
$a_2$ & $743646$ & $573$ & $2342536$ \\
$a_3$ & $820934$ & $722$ & $4569115$ \\
$a_4$ & $1000000$ & $1000$ & $9077428$ \\
$a_5$ & $613675$ & $695$ & $8382162$ \\
$a_6$ & $541754$ & $639$ & $10000000$ \\
$a_9$ & $-252836$ & $-220$ & $-5494550$ \\
$a_{10}$ & $-81234$ & $0$ & $0$ \\
$a_{11}$ & $-94750$ & $0$ & $0$ \\
\midrule
$M(d,3)$ & $30$ & $50$ & $275$ \\
$f(1)/a_0$ & $32.0003\ldots$ & $53.2$ & $282.0011\ldots$ \\
lower bound for $N$ & $33$ & $54$ & $283$ \\
\bottomrule
\end{tabular}
\end{table}

\end{appendices}


\begin{thebibliography}{9}
\bibitem{NN03}
S.~Nikova, V.~Nikov,
Improvement of the Delsarte bound for $\tau$-designs when it is not the best bound possible, {\it Des. Codes Cryptogr.} \textbf{28} (2003), 201--222.
\end{thebibliography}
\end{document}